\documentclass[11pt]{article}
\usepackage{latexsym,amssymb}
\usepackage{amsmath}
\usepackage[mathscr]{eucal}
\usepackage{color}
\usepackage[abbrev]{amsrefs}

\newtheorem{Theorem}{\bf Theorem}[section]
\newtheorem{Lemma}{\bf Lemma}[section]
\newtheorem{Proposition}{\bf Proposition}[section]  
\newtheorem{Corollary}{\bf Corollary}[section]
\newtheorem{Remark}{\bf Remark}[section]
\newtheorem{Example}{\bf Example}[section]
\newtheorem{Definition}{\bf Definition}[section]

\newenvironment{theorem}{\begin{Theorem}$\!\!\!$}{\end{Theorem}}
\newenvironment{lemma}{\begin{Lemma}$\!\!\!$}{\end{Lemma}}
\newenvironment{proposition}{\begin{Proposition}$\!\!\!$}{\end{Proposition}}

\newenvironment{remark}{\begin{Remark}$\!\!\!$}{\end{Remark}}

\numberwithin{equation}{section}

\def\XXint#1#2#3{{\setbox0=\hbox{$#1{#2#3}{\int}$}
\vcenter{\hbox{$#2#3$}}\kern-.5\wd0}}

\begin{document}
\title{Decay estimates for Schr\"odinger heat semigroup\\ 
with inverse square potential in Lorentz spaces II}
\author{Kazuhiro Ishige and Yujiro Tateishi\\
Graduate School of Mathematical Sciences, The University of Tokyo\\
3-8-1 Komaba, Meguro-ku, Tokyo 153-8914, Japan\vspace{10pt}
}
\date{}
\maketitle
\noindent
\qquad\quad
e-mail address:

\noindent
\qquad\quad
ishige@ms.u-tokyo.ac.jp (K. Ishige), 
tateishi@ms.u-tokyo.ac.jp (Y. Tateishi)
$$
\qquad
$$
\begin{abstract}
Let $H:=-\Delta+V$ be a nonnegative Schr\"odinger operator on $L^2({\bf R}^N)$, 
where $N\ge 2$ and $V$ is a radially symmetric inverse square potential. 
Let $\|\nabla^\alpha e^{-tH}\|_{(L^{p,\sigma}\to L^{q,\theta})}$ be the operator norm of $\nabla^\alpha e^{-tH}$ 
from the Lorentz space $L^{p,\sigma}({\bf R}^N)$ to $L^{q,\theta}({\bf R}^N)$, 
where $\alpha\in\{0,1,2,\dots\}$.
We establish both of upper and lower decay estimates of $\|\nabla^\alpha e^{-tH}\|_{(L^{p,\sigma}\to L^{q,\theta})}$ 
and study sharp decay estimates of $\|\nabla^\alpha e^{-tH}\|_{(L^{p,\sigma}\to L^{q,\theta})}$. 
Furthermore, we characterize the Laplace operator $-\Delta$ 
from the view point of the decay of $\|\nabla^\alpha e^{-tH}\|_{(L^{p,\sigma}\to L^{q,\theta})}$.
\end{abstract}
\newpage
\section{Introduction}
This paper is concerned with the decay of derivatives of 
Schr\"odinger heat semigroups $e^{-tH}$, 
where $H:=-\Delta+V$ is a nonnegative Schr\"odinger operator in $L^2({\bf R}^N)$. 
Throughout this paper we assume 
the following condition~($\mbox{V}_m$), where $m\in\{1,2,\dots,\infty\}$: 
\begin{equation}
\tag{$\mbox{V}_m$}
\left\{
\begin{array}{ll}
({\rm i}) & \mbox{$V=V(|x|)$ in ${\bf R}^N\setminus\{0\}$ and $V\in C^m((0,\infty))$};\vspace{7pt}\\
({\rm ii}) & \mbox{$V(r)=\lambda_1r^{-2}+O(r^{-2+\rho_1})$ as $r\to +0$},\vspace{3pt}\\
 & \mbox{$V(r)=\lambda_2r^{-2}+O(r^{-2-\rho_2})$ as $r\to \infty$},\vspace{3pt}\\
  & \mbox{for some $\lambda_1$, $\lambda_2\in[\lambda_*,\infty)$ with $\lambda_*:=-(N-2)^2/4$ and $\rho_1$, $\rho_2>0$};\vspace{7pt}\\
({\rm iii}) &  \mbox{$\displaystyle{\sup_{r>0}\,\left|\,r^{\ell+2}\frac{d^\ell}{dr^\ell}V(r)\right|<\infty}$ for $\ell\in\{1,\dots,m\}$},
\end{array}
\right.
\end{equation}
and investigate decay estimates of $\|\nabla^\alpha e^{-tH}\|_{(L^{p,\sigma}\to L^{q,\theta})}$, 
where $\alpha\in\{0,1,\dots,m+1\}$. 
Here $\|\nabla^\alpha e^{-tH}\|_{(L^{p,\sigma}\to L^{q,\theta})}$ is the operator norm 
from the Lorentz space $L^{p,\sigma}({\bf R}^N)$ to $L^{q,\theta}({\bf R}^N)$, that is, 
$$
\|\nabla^\alpha e^{-tH}\|_{(L^{p,\sigma}\to L^{q,\theta})}
:=\sup
\left\{\left\|\nabla^\alpha e^{-tH}\phi\right\|_{L^{q,\theta}({\bf R}^N)}\,:\,
\mbox{$\phi\in C_{\rm c}({\bf R}^N)$ with $\|\phi\|_{L^{p,\sigma}({\bf R}^N)}=1$}\right\}.
$$
Here
$$
(p,q,\sigma,\theta)\in\Lambda:=
\left\{1\le p\le q\le\infty,\,\sigma,\,\theta\in[1,\infty]:
\begin{array}{ll}
\mbox{$\sigma=1$\, if \,$p=1$}, & \mbox{$\sigma=\infty$\, if \,$p=\infty$}\vspace{3pt}\\
\mbox{$\theta=1$\, if \,$q=1$}, & \mbox{$\theta=\infty$\, if \,$q=\infty$}\vspace{3pt}\\
\mbox{$\sigma\le\theta$\, if \,$p=q$} &
\end{array}
\right\}.
$$

Nonnegative Schr\"odinger operators $H:=-\Delta+V$ on $L^2({\bf R}^N)$ and their heat semigroups $e^{-tH}$ 
have been studied by many mathematicians since the pioneering work due to Simon~\cite{S}. 
See e.g. \cite{Bar}, \cite{CK}, \cite{CUR}, \cite{DS}, 
\cite{IIY01}, \cite{IIY02}, \cite{IK01}--\cite{IKO}, 
\cite{LS}--\cite{M}, 
\cite{P3}--\cite{Zhang}, and references therein.   
(See also the monographs of Davies \cite{Dav}, Grigor'yan \cite{Gri}, and Ouhabaz \cite{Ouh}.) 
The Schr\"odinger operator with an inverse square potential often appears in the field of Schr\"odinger operators 
and nonlinear PDEs such as semilinear parabolic equations, 
and the decay estimates of $e^{-tH}\phi$ and their derivatives are crucial 
for the study of the behavior of $e^{-tH}\phi$. 
For related results on the decay of $\|\nabla^\alpha e^{-tH}\|_{(L^{p,\sigma}\to L^{q,\theta})}$, 
see e.g. \cite{LL}, \cite{DS}, \cite{Ishige}, \cite{IK01}, \cite{IK05}, \cite{IT01}, \cite{S}, and references therein.  

This paper is a continuation of our previous paper~\cite{IT01}, 
where the authors of this paper obtained upper decay estimates of $\|\nabla e^{-tH}\|_{(L^{p,\sigma}\to L^{q,\theta})}$, 
where $(p,q,\sigma,\theta)\in\Lambda$, under condition~($\mbox{V}_1$). 
In this paper we develop the arguments in~\cite{IT01} and obtain upper decay estimates of 
$\|\nabla^\alpha e^{-tH}\|_{(L^{p,\sigma}\to L^{q,\theta})}$ systematically, where $\alpha\in\{0,1,\dots,m+1\}$. 
Furthermore, we also establish lower decay estimates of $\|\nabla^\alpha e^{-tH}\|_{(L^{p,\sigma}\to L^{q,\theta})}$.
Combining both of upper and lower decay estimates, we study sharp decay estimates of 
$\|\nabla^\alpha e^{-tH}\|_{(L^{p,\sigma}\to L^{q,\theta})}$.
Furthermore, we give a new characterization of $-\Delta$ from the view point of the decay of 
$\|\nabla^\alpha e^{-tH}\|_{(L^{p,\sigma}\to L^{q,\theta})}$ (see Theorem~\ref{Theorem:1.3}).

We introduce some notations. 
Set ${\bf N}_0:=\{0,1,2,\dots\}$. 
Let $B(x,R):=\{y\in{\bf R}^N\,:\,|y-x|<R\}$ and $B(x,R)^c:={\bf R}^N\setminus B(x,R)$ for $x\in{\bf R}^N$ and $R>0$.
For any $r\in[1,\infty]$, let $r'$ be the H\"older conjugate number of $r$, that is, 
$$
r'=\frac{r}{r-1}\quad\mbox{if}\quad 1<r<\infty,
\quad
r'=1\quad\mbox{if}\quad r=\infty,
\quad
r'=\infty\quad\mbox{if}\quad r=1.
$$
Let $\Delta_{{\bf S}^{N-1}}$ be the Laplace-Beltrami operator on ${\bf S}^{N-1}$. 
Let $\{\omega_k\}_{k=0}^\infty$ be 
the eigenvalues of 
\begin{equation}
\label{eq:1.1}
-\Delta_{{\bf S}^{N-1}}Q=\omega Q\quad\mbox{on}\quad{\bf S}^{N-1},
\qquad
Q\in L^2({\bf S}^{N-1}). 
\end{equation}
Then $\omega_k=k(N+k-2)$ for $k=0,1,2,\dots$. 
Let
$\{Q_{k,i}\}_{i=1}^{d_k}$ and $d_k$ be 
the orthonormal system and the dimension of the eigenspace corresponding to $\omega_k$, respectively. 
Here 
\begin{equation}
\label{eq:1.2}
d_k=\frac{(N+2k-2)(N+k-3)!}{(N-2)!k!}=O(k^{N-2})\quad\mbox{as}\quad k\to\infty.
\end{equation}
Assume condition~($\mbox{V}_m$) 
and let $H:=-\Delta+V$ be nonnegative, that is, 
\begin{equation}
\tag{N}
\int_{{\bf R}^N}\left[|\nabla\phi|^2+V(|x|)\phi^2\right]\,dx\ge 0,
\qquad\phi\in C_c^\infty({\bf R}^N\setminus\{0\}). 
\end{equation}
The operator $H$ is said {\it subcritical} if, for any $W\in C_c({\bf R}^N)$, 
$H-\epsilon W$ is nonnegative for small enough $\epsilon>0$. 
If not, the operator $H$ is said {\it critical}. 
For any $k\in\{0,1,2,\dots\}$, set $A_{1,k}:=A^+_{\lambda_1+\omega_k}$ and 
\begin{equation}
\label{eq:1.3}
\begin{split}
A_{2,k}:= & 
\left\{
\begin{array}{ll}
A^-_{\lambda_2} & \mbox{if $k=0$ and $H$ is critical},\vspace{3pt}\\
A^+_{\lambda_2+\omega_k} & \mbox{otherwise},
\end{array}
\right.\\
B_k:= & 
\left\{
\begin{array}{ll}
1 &  \mbox{if $k=0$, $\lambda_2=\lambda_*$ and $H$ is subcritical},\vspace{3pt}\\
0 &  \mbox{otherwise}.
\end{array}
\right.
\end{split}
\end{equation}
Here 
\begin{equation}
\label{eq:1.4}
A^\pm_\lambda:=\frac{-(N-2)\pm\sqrt{D_\lambda}}{2}
\quad\mbox{for $\lambda\ge\lambda_*$, 
where $D_\lambda:=(N-2)^2+4\lambda$}.
\end{equation}
By the standard theory for ordinary differential equations 
we see that, 
for any $k\in\{0,1,2,\dots\}$, 
there exists a unique solution $h_k=h_k(r)$ to the problem
\begin{equation}
\label{eq:1.5}
\begin{split}
 & h_k''+\frac{N-1}{r}h_k'-V_k(r)h_k=0\quad\mbox{in}\quad(0,\infty),\\
 & h_k(r)=r^{A_{1,k}}(1+o(1))\quad\mbox{as}\quad r\to+0,
\end{split}
\end{equation}
where $V_k(r):=V(r)+\omega_kr^{-2}$. 
(See also Section~2.1.) 
Notice that $h_k\in L^2(B(0,1))$. 
Furthermore, it follows from \cite[Theorem~1.1]{IKO} that 
$h_k(r)>0$ for $r>0$ and 
\begin{equation}
\label{eq:1.6}
h_k(r)=c_kv_k(r)(1+o(1))\quad\mbox{as}\quad r\to\infty,
\quad\mbox{where}\quad
v_k(r):=r^{A_{2,k}}(\log r)^{B_k},
\end{equation}
for some $c_k>0$.
For $(k,i)\in{\mathcal K}:=\{(k,i)\,:\,k=0,1,2,\dots,\,\,i=1,\dots,d_k\}$, setting 
\begin{equation}
\label{eq:1.7}
J_{k,i}(x):=h_k(|x|)Q_{k,i}\left(\frac{x}{|x|}\right),
\end{equation}
we see that $J_{k,i}$ is a harmonic function for $H$, 
that is, $HJ_{k,i}=0$ in ${\bf R}^N\setminus\{0\}$. 
In particular,  
the function~$h_0$ is said {\it a positive harmonic function} for the operator $H$. 
When $H$ is critical, if $h_0\not\in L^2({\bf R}^N)$, 
then $H$ is said {\it null-critical}\,: if not, $H$ is said {\it positive-critical}. 
The decay of the fundamental solution $p=p(x,y,t)$ corresponding to $e^{-tH}$ 
depends on whether $H$ is either subcritical, null-critical or positive-critical. 
In particular, if $H$ is positive-critical, then $e^{-tH}\phi$ does not necessarily decay as $t\to\infty$.  
See \cite{P3}. (See also \cite{IKO}.)
\medskip

In this paper, under condition~($\mbox{V}_m$), 
we assume either
\begin{equation}
\tag{N'}
\qquad\qquad
{\rm (i)}\quad\mbox{$H$ is subcritical}\quad\mbox{or}\quad
{\rm (ii)}\quad\mbox{$H$ is critical and $A_{2,0}>-N/2$},
\qquad
\end{equation}
and obtain both of upper and lower decay estimates of $\|\nabla^\alpha e^{-tH}\|_{(L^{p,\sigma}\to L^{q,\theta})}$,  
where $\alpha\in\{0,1,\dots,m+1\}$ and $(p,q,\sigma,\theta)\in\Lambda$. 
Case~(ii) is in the null-critical one (see also \cite[Remark~1.1~(iii)]{IT01}). 

We state some results in this paper. 
These are obtained as applications of our upper and lower decay estimates in Sections~3 and 4. 
In what follows, for $(p,p,\sigma,\sigma)\in\Lambda$ and $t>0$, set 
$$
\Gamma_{p,\sigma}(t):=\frac{\|h_0\|_{L^{p,\sigma}(B(0,\sqrt{t}))}}{h_0(\sqrt{t})}\quad\mbox{if}\quad h_0\in L^{p,\sigma}(B(0,1)),
\quad
\Gamma_{p,\sigma}(t):=\infty\quad\mbox{if}\quad h_0\not\in L^{p,\sigma}(B(0,1)).
$$
The first theorem clarifies the relationship between the decay of 
$\|\nabla^\alpha e^{-tH}\|_{(L^{p,\sigma}\to L^{q,\theta})}$ and harmonic functions for $H$, 
where $\alpha\in\{0,1,2\}$ and $(p,q,\theta,\sigma)\in\Lambda$. 
\begin{theorem}
\label{Theorem:1.1}
Assume conditions~{\rm ($\mbox{V}_1$)} and {\rm (N')}.  
Let $(p,q,\sigma,\theta)\in\Lambda$ and $\alpha\in\{0,1,2\}$. 
Then there exists $C>0$ such that 
\begin{equation}
\label{eq:1.8}
C^{-1}\Phi_\alpha(t)\le\|\nabla^\alpha e^{-tH}\|_{(L^{p,\sigma}\to L^{q,\theta})}\le C\Phi_\alpha(t)
\quad\mbox{for}\quad t>0.
\end{equation}
Here
$$
\Phi_\alpha(t):=
\left\{
\begin{array}{ll}
t^{-\frac{N}{2}}\Gamma_{p',\sigma'}(t)\Gamma_{q,\theta}(t)
 & \mbox{if}\quad \alpha=0,\vspace{7pt}\\
\displaystyle{t^{-\frac{N}{2}}\Gamma_{p',\sigma'}(t)
\left[\frac{\|\nabla h_0\|_{L^{q,\theta}(B(0,\sqrt{t}))}}{h_0(\sqrt{t})}+t^{\frac{N}{2q}-\frac{1}{2}}\right]} 
 & \mbox{if}\quad \alpha=1,\vspace{7pt}\\
 \displaystyle{t^{-\frac{N}{2}}\Gamma_{p',\sigma'}(t)
\left[\frac{\|\nabla^2 h_0\|_{L^{q,\theta}(B(0,\sqrt{t}))}}{h_0(\sqrt{t})}+\sum_{i=1}^N\frac{\|\nabla^2 J_{1,i}\|_{L^{q,\theta}(B(0,\sqrt{t}))}}{h_1(\sqrt{t})}
+t^{\frac{N}{2q}-1}\right]} 
 & \mbox{if}\quad \alpha=2.
\end{array}
\right.
$$
\end{theorem}
\begin{remark}
\label{Remark:1.1}
Let $V\in C^1([0,\infty))$ and assume conditions~{\rm ($\mbox{V}_1$)} with $\lambda_1, \lambda_2\in(\lambda_*,\infty)$ and {\rm (N')}.  
Then the sharp large time decay estimate of $\|e^{-tH}\|_{(L^{p,\sigma}\to L^{q,\theta})}$ has been already obtained in \cite{IIY02}. 
Our decay estimate~\eqref{eq:1.8} with $\alpha=0$ gives the same decay estimate as in \cite{IIY02} 
and it has a simpler expression. 
\end{remark}
In the second and the third theorems, under conditions $(\mbox{V}_m$) and (N'), 
we characterize the Schr\"odinger operator $H$ satisfying 
$$
\|\nabla^\alpha e^{-tH}\|_{(L^{p,\sigma}\to L^{q,\theta})}
\le Ct^{-\frac{N}{2}\left(\frac{1}{p}-\frac{1}{q}\right)-\frac{\alpha}{2}},
\quad t>0,
$$
for some $C>0$, where $\alpha\in\{0,1,\dots,m\}$ and $(p,q,\sigma,\theta)\in\Lambda$.  
\begin{theorem}
\label{Theorem:1.2}
Let $m\in\{1,2,\dots\}$ and assume conditions~{\rm ($\mbox{V}_m$)} and {\rm (N')}. 
\begin{itemize}
  \item[{\rm (a)}] 
  For any $(p,q,\sigma,\theta)\in\Lambda$ and $\alpha\in\{0,1,\dots,m+1\}$, 
  there exists $C_1>0$ such that 
  $$
  \|\nabla^\alpha e^{-tH}\|_{(L^{p,\sigma}\to L^{q,\theta})}
  \ge C_1^{-1}t^{-\frac{N}{2}\left(\frac{1}{p}-\frac{1}{q}\right)-\frac{\alpha}{2}},\qquad t>0.
  $$
  \item[{\rm (b)}]
  Let $V\not\equiv 0$ in $(0,\infty)$ and $\alpha\in\{0,1,\dots,m+1\}$. 
  Assume that
  there exists $C_2>0$ such that 
  \begin{equation}
  \label{eq:1.9}
  \|\nabla^\alpha e^{-tH}\|_{(L^{p,\sigma}\to L^\infty)}
  \le C_2t^{-\frac{N}{2p}-\frac{\alpha}{2}},\qquad t>0,
  \end{equation}
  for some $(p,\infty,\sigma,\infty)\in\Lambda$. 
  Then $H$ must be subcritical with $\lambda_1\in[\omega_\alpha,\infty)\cup\{0\}$ and $\lambda_2\in[\omega_\alpha,\infty)$.
\end{itemize}
\end{theorem}
\begin{theorem}
\label{Theorem:1.3}
Assume conditions~{\rm ($\mbox{V}_\infty$)} and {\rm (N')}. 
Let $(p,q,\sigma,\theta)\in\Lambda$. 
Assume that, for any $\alpha\in\{0,1,2,\dots\}$, 
there exists $C>0$ such that 
\begin{equation}
\label{eq:1.10}
\|\nabla^\alpha e^{-tH}\|_{(L^{p,\sigma}\to L^{q,\theta})}
\le Ct^{-\frac{N}{2}\left(\frac{1}{p}-\frac{1}{q}\right)-\frac{\alpha}{2}},\quad t>0. 
\end{equation}
Then $V$ must be identically zero in ${\bf R}^N$, that is, $H=-\Delta$. 
\end{theorem}
See also Theorems~\ref{Theorem:7.3} and \ref{Theorem:7.4}.
\vspace{3pt}

Upper and lower decay estimates of $\|\nabla^\alpha e^{-tH}\|_{(L^{p,\sigma}\to L^{q,\theta})}$ 
are given in Sections~3 and 4, respectively. 
These are main ingredients of this paper. 
In order to obtain the upper decay estimates, 
we follow the arguments in \cite{IK01} and \cite{IT01}.
For any $\phi\in C_c({\bf R}^N)$, 
we find radially symmetric functions 
$\{\phi_{k,i}\}_{(k,i)\in{\mathcal K}}\subset L^2({\bf R}^N)$ 
such that 
\begin{equation}
\label{eq:1.11}
\phi(x)=\sum_{k=0}^\infty\sum_{i=1}^{d_k}\phi_{k,i}(|x|)Q_{k,i}\left(\frac{x}{|x|}\right)
\quad\mbox{in}\quad L^2({\bf R}^N).
\end{equation}
Let $H_k:=-\Delta+V_k(|x|)$ and set 
\begin{equation}
\label{eq:1.12}
v_{k,i}(|x|,t):=[e^{-tH_k}\phi_{k,i}](|x|),
\qquad
u_{k,i}(x,t):=v_{k,i}(|x|,t)Q_{k,i}\left(\frac{x}{|x|}\right).
\end{equation}
Then 
\begin{equation}
\label{eq:1.13}
\left[e^{-tH}\phi\right](x)=\sum_{k=0}^\infty\sum_{i=1}^{d_k} u_{k,i}(x,t)
=\sum_{k=0}^\infty\sum_{i=1}^{d_k} v_{k,i}(|x|,t)Q_{k,i}\left(\frac{x}{|x|}\right)
\quad\mbox{in}\quad C^2(K)
\end{equation}
for compact sets $K\subset{\bf R}^N\setminus\{0\}$ and $t>0$ 
(see \cite{IKM} and \cite{IM}). 
We study the behavior of derivatives of $v_{k,i}$ by using the radially symmetry of $v_{k,i}$, 
and obtain upper decay estimates of $\|\nabla^\alpha e^{-tH}\|_{(L^{p,\sigma}\to L^{q,\theta})}$. 
The introduction of functions $I_k^n[\cdot]$ and $J_{k,i}^n$ enables us to obtain 
those decay estimates systematically. 
See Subsection~2.2 and Section~3. 

On the other hand, 
we study lower decay estimates of $\|\nabla^\alpha e^{-tH}\|_{(L^{p,\sigma}\to L^{q,\theta})}$
by using the lower Gaussian estimate of the fundamental solution to a parabolic equation with $A_2$-weight 
and by applying a tricky switch of weight (see \eqref{eq:4.6}). 
See Section~4.

The rest of this paper is organized as follows. 
In Section~2 we recall some properties of Lorentz spaces and $h_k$.  
Furthermore, we obtain some preliminary results on $I_k^n[\cdot]$, $J_{k,i}^n$, and $e^{-tH}$. 
In Sections~3 and 4, 
we obtain upper and lower decay estimates of $\|\nabla^\alpha e^{-tH}\|_{(L^{p,\sigma}\to L^{q,\theta})}$, respectively. 
Section~5 is devoted to the proof of Theorem~\ref{Theorem:1.1}.  
In Section~6 we prove Theorems~\ref{Theorem:1.2} and \ref{Theorem:1.3}. 
In Section~7, as typical examples of inverse square potentials, 
we treat the Hardy potentials and bounded potentials, 
and clarify the relationship between 
the corresponding harmonic functions and the decay of $\|\nabla^\alpha e^{-tH}\|_{(L^{p,\sigma}\to L^{q,\theta})}$. 
Furthermore, we show that 
the decay of $\|\nabla^\alpha e^{-tH}\|_{(L^{p,\sigma}\to L^{q,\theta})}$ is delicate in the case when $\lambda_1=\lambda_2=0$.
\section{Preliminaries}
In this section we introduce Lorentz spaces. 
Furthermore, we recall some results on $h_k$ and $e^{-tH}$,
and prove some preliminary results. 

Throughout this paper we use the same definition of $e^{-tH}$ 
and notations as in \cite{IT01}. 
In particular, 
for any positive functions $f$ and $g$ on a set $E$, 
we write $f\asymp g$ for $x\in E$ if there exists $c>0$ such that 
$c^{-1}\le f(x)/g(x)\le c$ for $x\in E$. 
By the letters $C$ and $D$ we denote generic positive constants 
and they may have different values also within the same line. 
\subsection{Lorentz spaces}
For any measurable function $\phi$ in ${\bf R}^N$, 
we denote by $\mu=\mu(\lambda)$ the distribution function of $\phi$, that is, 
$$
\mu(\lambda):=\left|\{x\in{\bf R}^N\,:\,|\phi(x)|>\lambda\}\right| 
\quad\mbox{for}\quad \lambda > 0.
$$
Here $|E|$ is the $N$-dimensional Lebesgue measure of $E$ for measurable sets $E$ in ${\bf R}^N$.
We define 
the non-increasing rearrangement $\phi^*$ of $\phi$
and the spherical rearrangement $\phi^{\sharp}$ of $\phi$ by 
$$
\phi^{*}(s):=\inf\{\lambda>0\,:\,\mu(\lambda)\le s\},
\qquad
\phi^\sharp(x):=\phi^*(\alpha_N|x|^N),
$$
for $s>0$ and $x\in {\bf R}^N$, respectively,
where $\alpha_N$ is the volume of the unit ball in ${\bf R}^N$. 
For any $(p,p,\sigma,\sigma)\in\Lambda$, 
we define the Lorentz space $L^{p,\sigma}({\bf R}^N)$ by 
$$
L^{p,\sigma}({\bf R}^N):=\{\phi\,:\, \mbox{$\phi$ is measurable in ${\bf R}^N$},\,\,\, \|\phi\|_{L^{p,\sigma}}<\infty\},
$$
where 
$$
\|\phi\|_{L^{p,\sigma}}:=
 \left\{
\begin{array}{ll}
  \displaystyle{\biggr(
   \int_{{\bf R}^N}\left(|x|^{\frac{N}{p}}\phi^{\sharp}(x)\right)^{\sigma}\frac{dx}{|x|^N}
  \biggr)^{\frac{1}{\sigma}}}\quad & \mbox{if}\quad 1\le \sigma<\infty, \vspace{5pt}\\
 \displaystyle{\sup_{x\in {\bf R}^N}\,|x|^{\frac{N}{p}}\phi^{\sharp}(x)}\quad & \mbox{if}\quad\sigma=\infty.
\end{array}
 \right.
$$
Here $N/p=0$ if $p=\infty$. 
The Lorentz spaces have the following properties:  
\begin{equation*}
\begin{array}{ll}
L^{p,p}({\bf R}^N)=L^p({\bf R}^N) & \mbox{if $1\le p\le \infty$};\vspace{3pt}\\ 
L^{p,\sigma_1}({\bf R}^N)\subset L^{p,\sigma_2}({\bf R}^N) & \mbox{if $1\le p<\infty$ and $1\le \sigma_1\le \sigma_2 \le \infty$}.\vspace{3pt}\\ 
\end{array}
\end{equation*}
Furthermore, 
there exists $C>0$ depending only on $N$ such that 
\begin{eqnarray*}
 & & \|f+g\|_{L^{p,\sigma}}\le C(\|f\|_{L^{p.\sigma}}+\|g\|_{L^{p,\sigma}})\quad\mbox{if $f,g\in L^{p,\sigma}({\bf R}^N)$},\\
 & & \|fg\|_{L^1}\le C\|f\|_{L^{p,\sigma}}\|g\|_{L^{p',\sigma'}}
 \qquad\qquad\,\,\mbox{if $f\in L^{p,\sigma}({\bf R}^N)$, $g\in L^{p',\sigma'}({\bf R}^N)$},\\
 & & \|f*g\|_{L^{q,\theta}}\le C\|f\|_{L^{p,\sigma}}\|g\|_{L^{r,s}}
 \qquad\quad\,\,\,\mbox{if $f\in L^{p,\sigma}({\bf R}^N)$, $g\in L^{r,s}({\bf R}^N)$}.
\end{eqnarray*}
Here $(q,q,\theta,\theta)$, $(r,r,s,s)\in\Lambda$ and
\begin{equation*}
   \frac{1}{r}+\frac{1}{p}=\frac{1}{q}+1,\qquad \frac{1}{\theta}=\frac{1}{s}+\frac{1}{\sigma}.
\end{equation*}
(See e.g. \cite{BS} and \cite{Gra}.) 
For any measurable function~$f$ in a domain $\Omega$, 
we say that $f\in L^{p,\sigma}(\Omega)$ if and only if $\tilde{f}\in L^{p,\sigma}({\bf R}^N)$, 
where $\tilde{f}$ is the zero extension of $f$ to ${\bf R}^N$. 
Furthermore, we write $\|f\|_{L^{p,\sigma}(\Omega)}=\|\tilde{f}\|_{L^{p,\sigma}}$.
Then, for any $A\in{\bf R}$, 
the function $f_A$ defined by $f_A(x):=|x|^A$ satisfies $f_A\in L^{p,\sigma}(B(0,R))$ for $R>0$
if and only if 
\begin{equation}
\label{eq:2.1}
 pA+N>0\quad\mbox{for}\quad 1\le\sigma<\infty,\qquad
 pA+N\ge 0\quad\mbox{for}\quad \sigma=\infty.
 \end{equation}
 Furthermore, for any $R>0$, under condition~\eqref{eq:2.1}, we have 
$$
\|f_A\|_{L^{p,\sigma}(B(0,\sqrt{t}))}\asymp t^{\frac{A}{2}+\frac{N}{2p}}\quad\mbox{for}\quad t\in(0,R^2].
$$
In particular, 
for any $k\in\{0,1,2,\dots\}$, 
by \eqref{eq:1.6} we see that 
$$
\frac{\|h_k\|_{L^{p,\sigma}(B(0,\sqrt{t}))}}{h_k(\sqrt{t})}\asymp t^{\frac{N}{2p}}\quad\mbox{for}\quad 0<t\le R^2\quad\mbox{if}\quad h_k\in L^{p.\sigma}(B(0,1)).
$$
\subsection{Preliminary results on $h_k$}
Consider the ordinary differential equation
\begin{equation}
\label{eq:2.2}
h''+\frac{N-1}{r}h'-V_k(r)h=0\quad\mbox{in}\quad(0,\infty).
\end{equation}
Then ODE~\eqref{eq:2.2} has two linearly independent solutions $h_k^+$ and $h_k^-$ such that 
$$
h_k^+(r)=v_{k,\lambda_1}^+(r)(1+o(1)),\qquad
h_k^-(r)=v_{k,\lambda_1}^-(r)(1+o(1)),
$$
as $r\to +0$ and $h_k^-(1)=1$.
Here
$$
v_{k,\lambda}^+(r):=r^{A^+_{\lambda+\omega_k}},
\quad 
v_{k,\lambda}^-(r):=\left\{
\begin{array}{ll}
r^{-\frac{N-2}{2}}\displaystyle{\left|\log\frac{r}{2}\right|} & \mbox{if $\lambda=\lambda_*$ and $k=0$},\vspace{5pt}\\
r^{A^-_{\lambda+\omega_k}}  & \mbox{otherwise},
\end{array}
\right.
$$
for $\lambda\ge\lambda_*$, where $A^\pm_\lambda$ is as in \eqref{eq:1.4}. 
Then we have the following two propositions for the solution~$h_k$ to \eqref{eq:1.5}. 
\begin{proposition}
\label{Proposition:2.1}
Let $m\in\{1,2,\dots\}$ and assume conditions~{\rm ($\mbox{V}_m$)} and {\rm (N')}. 
Then, for any $\ell\in\{0,1,\dots,m+1\}$, 
$$
\frac{d^\ell}{dr^\ell}h_k(r)=
\left\{
\begin{array}{ll}
\displaystyle{\frac{d^\ell}{dr^\ell}}v^+_{k,\lambda_1}(r)+O\left(r^{-\ell+\rho_1}v^+_{k,\lambda_1}(r)\right)
 & \quad\mbox{as}\quad r\to +0,\vspace{5pt}\\
c_k\displaystyle{\frac{d^\ell}{dr^\ell}}v_{k}(r)+o\left(r^{-\ell}v_{k}(r)\right)
 & \quad\mbox{as}\quad r\to\infty.
\end{array}
\right.
$$
Here $v_k$ is as in \eqref{eq:1.6}.
Furthermore, there exists $C>0$ such that 
\begin{equation*}
\begin{split}
 & C^{-1}\le\displaystyle{\frac{h_k(r)}{v_{k,\lambda_1}^+(r)}}\le C
\quad\mbox{in}\quad(0,1],\qquad
C^{-1}\le\displaystyle{\frac{h_k(r)}{v_k(r)}}\le C\quad\mbox{in}\quad(1,\infty),\\
 & \left|\frac{d^\ell}{dr^\ell}h_k(r)\right|\le C(k+1)^{\ell-1}r^{-\ell}h_k(r)\quad\mbox{in}\quad(0,\infty),
\end{split}
\end{equation*}
for $k\in\{0,1,2,\dots\}$ and $\ell\in\{0,1,\dots,m+1\}$. 
\end{proposition}
{\bf Proof.}
Proposition~\ref{Proposition:2.1} with $\ell\in\{0,1\}$ follows from \cite[Propositions~2.1 and 2.2]{IT01}. 
Then, by the use of equation \eqref{eq:2.2} we obtain the other desired relations.  
$\Box$
\begin{proposition}
\label{Proposition:2.2} 
Assume conditions~{\rm ($\mbox{V}_1$)} and {\rm (N')}.
\begin{itemize}
  \item[{\rm (a)}] 
  Let $(p,p,\sigma,\sigma)\in\Lambda$ be such that $h_0\in L^{p,\sigma}(B(0,1))$. 
  There exists $C_1>0$, independent of $(p,\sigma)$, such that 
  $$
  \Gamma_{p,\sigma}(t)\ge C_1t^{\frac{N}{2p}}\quad\mbox{for}\quad t>0.
  $$ 
  \item[{\rm (b)}] 
  There exists $C_2>0$ such that 
  $$
  \int_0^r s^{N-1}h_k(s)^2\,ds\le C_2(k+1)^{-1}r^N h_k(r)^2
  $$
  for $r>0$ and $k\in\{0,1,2,\dots\}$. 
  \item[{\rm (c)}]
  Let $\ell\in\{0,1,2,\dots\}$. 
  Then there exist $C_3>0$ and $\gamma>0$ such that 
  $$
  \frac{h_k(\epsilon r)}{h_\ell(\epsilon r)}\le C\epsilon^{\left(\frac{k}{2}-\gamma\right)_+} \frac{h_k(r)}{h_\ell(r)}
  $$
  for $r>0$, $\epsilon\in(0,1)$ and $k\in\{\ell+1,\ell+2,\dots\}$. 
  Here $a_+:=\max\{a,0\}$ for $a\in{\bf R}$.
\end{itemize}
\end{proposition}
{\bf Proof.} 
Assertions~(a) and (b) follow from \cite[Proposition 2.4]{IT01}. 
It suffices to prove assertion~(c). 

We consider the case of $\ell\in\{1,2,\dots\}$.   
Set
$$
\iota_{\ell,k}(r):=\left\{
\begin{array}{ll}
r^{A_{1,k}-A_{1,\ell}} & \mbox{for}\quad \mbox{$0<r<1$},\vspace{5pt}\\
r^{A_{2,k}-A_{2,\ell}} & \mbox{for}\quad \mbox{$r\ge 1$},
\end{array}
\right.
$$
where $k\in\{\ell+1,\ell+2,\dots\}$. 
Since $A_{i,k}\ge A_{i,\ell}$ for $i\in\{1,2\}$, 
the function $\iota_{\ell,k}$ is monotone increasing in $(0,\infty)$ for $k>\ell$. 
Furthermore, Proposition~\ref{Proposition:2.1} implies that
$$
\frac{h_k(r)}{h_\ell(r)}\asymp\iota_{\ell,k}(r)\quad\mbox{for $r>0$ and $k\in\{\ell+1,\ell+2,\dots\}$}. 
$$ 
Since $A_{i,k}=k(1+o(1))$ as $k\to\infty$, where $i=1,2$, 
we find $\gamma>0$ such that 
$$
A_{i,k}-A_{i,\ell}\ge\left(\frac{k}{2}-\gamma\right)_+\quad\mbox{for $k\in\{\ell+1,\ell+2,\dots\}$ and $i\in\{1,2\}$},
$$
which implies that 
$$
\frac{h_k(\epsilon r)}{h_\ell(\epsilon r)}\asymp 
\iota_{\ell,k}(\epsilon r)\le C\epsilon^{\left(\frac{k}{2}-\gamma\right)_+}\iota_{\ell,k}(r)
\asymp \epsilon^{\left(\frac{k}{2}-\gamma\right)_+}\frac{h_k(r)}{h_\ell(r)}
$$
for $r>0$, $\epsilon\in(0,1)$, and $k\in\{\ell+1,\ell+2,\dots\}$. 
Indeed, the inequality in the above relation holds since
\begin{equation*}
\begin{split}
 & \iota_{\ell,k}(\epsilon r)\le \epsilon^{A_{1,k}-A_{1,\ell}}\iota_{\ell,k}(r)\le\epsilon^{\left(\frac{k}{2}-\gamma\right)_+}\iota_{\ell,k}(r)
 \quad\mbox{if}\quad r\le 1,\\
 & \iota_{\ell,k}(\epsilon r)\le \epsilon^{A_{2,k}-A_{2,\ell}}\iota_{\ell,k}(r)\le\epsilon^{\left(\frac{k}{2}-\gamma\right)_+}\iota_{\ell,k}(r)
  \quad\mbox{if}\quad \epsilon r\ge 1,\\
 & \iota_{\ell,k}(\epsilon r)\le \epsilon^{A_{1,k}-A_{1,\ell}}\iota_{\ell,k}(1)
\le \epsilon^{\left(\frac{k}{2}-\gamma\right)_+}\iota_{\ell,k}(r)
\quad\mbox{if}\quad \epsilon r\le 1\le r.\\
\end{split}
\end{equation*}
Thus assertion~(c) follows for the case of $\ell\in\{1,2,\dots\}$. 

Let $\ell=0$. Similarly to the above argument, 
by \eqref{eq:1.5} and \eqref{eq:1.6} we see such that 
\begin{equation}
\label{eq:2.3}
\frac{h_1(\epsilon r)}{h_0(\epsilon r)}\le C\frac{h_1(r)}{h_0(r)}
\end{equation}
for $r>0$ and $\epsilon\in(0,1)$. 
Then, by assertion~(c) with $\ell\in\{1,2,\dots\}$ and \eqref{eq:2.3}, 
we find $\gamma'>0$ such that 
$$
\frac{h_k(\epsilon r)}{h_0(\epsilon r)}=\frac{h_k(\epsilon r)}{h_1(\epsilon r)}\frac{h_1(\epsilon r)}{h_0(\epsilon r)}
\le C\epsilon^{\left(\frac{k}{2}-\gamma'\right)_+}\frac{h_k(r)}{h_1(r)}\frac{h_1(r)}{h_0(r)}
=C\epsilon^{\left(\frac{k}{2}-\gamma'\right)_+}\frac{h_k(r)}{h_0(r)}
$$
for $r>0$, $\epsilon\in(0,1)$, and $k\in\{1,2,\dots\}$. 
This implies assertion~(c) for the case of $\ell=0$. 
Thus Proposition~\ref{Proposition:2.2} follows.
$\Box$
\vspace{5pt}

In the study of the decay of $\|\nabla^\alpha e^{-tH}\|_{(L^{p,\sigma}\to L^{q,\theta})}$, 
it is crucial to obtain the estimates of the spatial derivatives of $e^{-tH}\phi$ systematically. 
For this aim, we introduce functions $I_k^n[\cdot]$ and $J_{k,i}^n$. 
For any continuous function $f$ in $(0,\infty)$, 
we set 
$$
I_k[f](r):=\int_0^r s^{-N+1}h_k(s)^{-2}\left(\int_0^s \tau^{N-1}h_k(\tau)^2f(\tau)\,d\tau\right)\,ds,
\quad k\in\{0,1,2,\dots\}.
$$
Then, under a suitable assumption of $f$ at $r=0$, 
we have
\begin{equation*}
\begin{split}
 & \frac{d^2}{dr^2}I_k[f](r)+\frac{N-1}{r}\frac{d}{dr}I_k[f](r)-\left(V(r)+\frac{\omega_k}{r^2}\right)I_k[f](r)=f(r)
\quad\mbox{in}\quad(0,\infty),\\
 & I_k[f](0)=\frac{d}{dr}I_k[f](0)=0.
\end{split}
\end{equation*}
We define $I_k^n[f]$, where $n\in\{0,1,2,\dots\}$, by
$$
I_k^{n+1}[f](r):=I_k[I_k^n[f]](r),\qquad I_k^0[f](r):=f(r).
$$
Furthermore, we set 
$$
I_k^n(|x|):=I_k^n[1](|x|),\quad 
J_{k,i}^n(x):=h_k(|x|)I_k^n(|x|)Q_{k,i}\left(\frac{x}{|x|}\right),
\quad\mbox{for}\quad x\in{\bf R}^N\setminus\{0\}.
$$
Notice that $J_{k,i}(x)=J_{k,i}^0(x)$ (see \eqref{eq:1.7}). 
Then we have: 
\begin{lemma}
\label{Lemma:2.1}
Assume the same conditions as in Proposition~{\rm\ref{Proposition:2.1}}. 
Let $\ell\in\{0,1,\dots,m+1\}$ and $n\in\{0,1,2,\dots\}$.
\begin{itemize}
  \item[{\rm (a)}] 
  If $\ell\le 2n$, then there exist $C_1>0$ and $D_1>0$ such that 
  $$
  \left|\nabla^\ell I_k^n[f](|x|)\right|\le C_1(k+1)^{D_1}|x|^{2n-\ell}\sup_{0<s<|x|}|f(s)|,
  \quad x\in{\bf R}^N\setminus\{0\},
  $$
  for $f\in C([0,\infty))$.
  \item[{\rm (b)}]
  There exist positive constants $C_2$, $C_3$, and $D_2$ such that 
  \begin{equation}
  \label{eq:2.4}
  t^{-n}\frac{|\nabla^\ell J_{k,i}^n(x)|}{h_k(\sqrt{t})}\le C_2(k+1)^{D_2}t^{-n}|x|^{2n-\ell} \frac{h_k(|x|)}{h_k(\sqrt{t})}
  \le C_3(k+1)^{D_2}|x|^{-\ell} \frac{h_0(|x|)}{h_0(\sqrt{t})}
  \end{equation}
  for $x\in B(0,\sqrt{t})$, $t>0$, and $(k,i)\in{\mathcal K}$. 
\end{itemize}
\end{lemma}
{\bf Proof.}
Applying elliptic regularity theorems to \eqref{eq:1.1}, 
for any $\beta\in\{0,1,\dots,\ell\}$, 
we have 
\begin{equation}
\label{eq:2.5}
\left|\nabla^\beta Q_{k,i}\left(\frac{x}{|x|}\right)\right|\le C(k+1)^D |x|^{-\beta},
\quad x\in{\bf R}^N\setminus\{0\},
\,\,\,(k,i)\in{\mathcal K}. 
\end{equation}
Furthermore, it follows from Proposition~\ref{Proposition:2.2}~(c) that 
$$
\frac{h_k(|x|)}{h_k(\sqrt{t})}\le C\frac{h_0(|x|)}{h_0(\sqrt{t})},
\quad x\in{\bf R}^N\setminus\{0\}.
$$
Then Proposition~\ref{Proposition:2.1} together with Proposition~\ref{Proposition:2.2}~(b) 
implies the desired inequalities. 
$\Box$
\subsection{Estimates of $J_{k,i}^n$}
We collect estimates of $J_{k,i}^n$, which are used in the rest of this paper. 
Due to \eqref{eq:1.5}, for $\alpha\in\{0,1,\dots,m+1\}$, 
we divide the behavior of $h_0$ near $0$ into the following two cases, 
which depends on whether $\partial_r^\alpha h_0$ is degenerate at $r=0$ or not:
\begin{itemize}
  \item[(${\mathcal A}_\alpha$)] $h_0(r)\not\asymp r^A$ as $r\to 0$ for $A\in\{0,1,\dots\}$ with $A\le\alpha-1$;
  \item[(${\mathcal B}_\alpha$)] $h_0(r)\asymp r^A$ as $r\to 0$, where $A\in\{0,1,\dots\}$ with $A\le\alpha-1$.
\end{itemize}
Let $(p,q,\sigma,\theta)\in\Lambda$, $(k,i)\in{\mathcal K}$, and $n\in\{0,1,2,\dots\}$. 
Set $h_k^{\langle\gamma\rangle}(r):=r^{-\gamma}h_k(r)$ for $\gamma\in{\bf R}$ and 
$$
\Gamma^k_{p,\sigma}(t):=\frac{\|h_k\|_{L^{p,\sigma}(B(0,\sqrt{t}))}}{h_k(\sqrt{t})}\quad\mbox{if}\quad h_k\in L^{p,\sigma}(B(0,1)),
\quad
\Gamma^k_{p,\sigma}(t):=\infty\quad\mbox{if}\quad h_k\not\in L^{p,\sigma}(B(0,1)).
$$
Then, by Proposition~\ref{Proposition:2.1}, Proposition~\ref{Proposition:2.2}, and Lemma~\ref{Lemma:2.1} 
we have: 
\vspace{5pt}
\newline
\underline{Case (${\mathcal A}_\alpha$)}: 
Consider case~(${\mathcal A}_\alpha$). 
By Proposition~\ref{Proposition:2.1} we find $R_1>0$ such that 
\begin{equation}
\label{eq:2.6}
|\nabla^\ell h_0|\le Ch_0^{\langle\ell\rangle}
\asymp \partial_r^\ell h_0\le |\nabla^\ell h_0| ,\qquad x\in B(0,R_1),
\end{equation}
for $\ell\in\{0,\dots,\alpha\}$. This together with \eqref{eq:2.4} implies that  
\begin{equation}
\label{eq:2.7}
t^{-n}\frac{\|\nabla^\ell J_{k,i}^n\|_{L^{q,\theta}(E)}}{h_k(\sqrt{t})}
\le C\frac{\|h_0^{\langle\ell\rangle}\|_{L^{q,\theta}(E)}}{h_0(\sqrt{t})}
\le C\frac{\|\nabla^\ell h_0\|_{L^{q,\theta}(E)}}{h_0(\sqrt{t})}
\end{equation}
for measurable sets $E\subset B(0,\sqrt{t})\cap B(0,R_1)$ and $\ell\in\{0,\dots,\alpha\}$.
\vspace{5pt}
\newline
\underline{Case (${\mathcal B}_\alpha$)}: 
Consider case~(${\mathcal B}_\alpha$). 
Let $R_2>0$. 
Then
\begin{equation}
\label{eq:2.8}
\Gamma_{p',\sigma'}(t)\asymp \Gamma^k_{p',\sigma'}(t)\asymp t^{\frac{N}{2p'}},
\qquad
\Gamma_{q,\theta}(t)\asymp \Gamma^k_{q,\theta}(t)\asymp t^{\frac{N}{2q}},
\end{equation}
for $0<t\le \sqrt{R_2}$. Furthermore, by \eqref{eq:2.4} we have
\begin{equation}
\label{eq:2.9}
t^{-n}\frac{\|\nabla^\ell J_{k,i}^n\|_{L^{q,\theta}(E)}}{h_k(\sqrt{t})}
\le Ct^{-\frac{\ell}{2}}|E|^{\frac{N}{q}}\le Ct^{\frac{N}{2q}-\frac{\ell}{2}}
\end{equation}
for measurable sets $E\subset B(0,\sqrt{t})\cap B(0,R_2)$ and $\ell\in\{0,\dots,\alpha\}$.
\vspace{5pt}

Similarly, we divide the behavior of $h_0$ at the space infinity into the following two cases: 
\begin{itemize}
  \item[(${\mathcal A}_\alpha'$)] $h_0(r)\not\asymp r^A$ as $r\to \infty$ for $A\in\{0,1,\dots\}$ with $A\le\alpha-1$;
  \item[(${\mathcal B}_\alpha'$)] $h_0(r)\asymp r^A$ as $r\to \infty$, where $A\in\{0,1,\dots\}$ with $A\le\alpha-1$.
\end{itemize}
\underline{Case (${\mathcal A}_\alpha'$)}:
Consider case~(${\mathcal A}_\alpha'$). 
By Proposition~\ref{Proposition:2.1} we find $R_3>0$ such that 
\begin{equation}
\label{eq:2.10}
|\nabla^\ell h_0|\le Ch_0^{\langle\ell\rangle}
\asymp \partial_r^\ell h_0\le |\nabla^\ell h_0| ,\qquad x\in B(0,R_3)^c,
\end{equation}
for $\ell\in\{0,\dots,\alpha\}$.
Then, similarly to \eqref{eq:2.7}, 
we have
\begin{equation}
\label{eq:2.11}
t^{-n}\frac{\|\nabla^\ell J_{k,i}^n\|_{L^{q,\theta}(E)}}{h_k(\sqrt{t})}
\le C\frac{\|h_0^{\langle\ell\rangle}\|_{L^{q,\theta}(E)}}{h_0(\sqrt{t})}
\le C\frac{\|\nabla^\ell h_0\|_{L^{q,\theta}(E)}}{h_0(\sqrt{t})}
\end{equation}
for measurable sets $E\subset B(0,\sqrt{t})\setminus B(0,R_3)$ and $\ell\in\{0,\dots,\alpha\}$.
\vspace{5pt}
\newline
\underline{Case (${\mathcal B}_\alpha'$)}: 
Consider case~(${\mathcal B}_\alpha'$). 
Let $R_4>0$. 
Then
\begin{equation}
\label{eq:2.12}
\Gamma_{p',\sigma'}(t)\asymp \Gamma^k_{p',\sigma'}(t)\asymp t^{\frac{N}{2p'}},
\quad
\Gamma_{q,\theta}(t)\asymp \Gamma^k_{q,\theta}(t)\asymp t^{\frac{N}{2q}},
\quad
h_k(\sqrt{t})\ge Ct^{\frac{k}{2}},
\end{equation}
for $t\ge \sqrt{R_4}$. Furthermore, by \eqref{eq:2.4} we have 
\begin{equation}
\label{eq:2.13}
t^{-n}\frac{\|\nabla^\ell J_{k,i}^n\|_{L^{q,\theta}(E)}}{h_k(\sqrt{t})}
\le Ct^{-\frac{\ell}{2}}|E|^{\frac{N}{q}}\le Ct^{\frac{N}{2q}-\frac{\ell}{2}}
\end{equation}
for measurable sets $E\subset B(0,\sqrt{t})\setminus B(0,R_4)$ and $\ell\in\{0,\dots,\alpha\}$.
\subsection{Estimates of solutions}
We obtain some results on the behavior of $e^{-tH}\phi$, where $\phi\in C_c({\bf R}^N)$. 
The following two propositions 
follow from the same arguments as in the proof of \cite[Propositions~3.1 and 4.2]{IT01}, respectively.
\begin{proposition}
\label{Proposition:2.3}
Let $m\in\{1,2,\dots\}$ and assume conditions~{\rm ($\mbox{V}_m$)} and {\rm (N')}.  
\begin{itemize}
  \item[{\rm (a)}] 
  Let $(p,q,\sigma,\theta)\in\Lambda$, $\alpha\in\{0,1,\dots,m+1\}$, $\beta\in\{0,1,2,\dots\}$, and $\delta\in (0,1]$.  
  Then there exist $C_1>0$ and $C_2>0$ such that 
  \begin{equation*}
  \begin{split}
   & t^\beta\left\|\partial_t^\beta\nabla^\alpha e^{-tH}\phi\right\|_{L^{q,\theta}(B(0,\delta\sqrt{t})^c)}\\
   & \le C_1t^{-\frac{N}{2}\left(1-\frac{1}{q}\right)-\frac{\alpha}{2}}
  \left[\frac{\|h_0\phi\|_{L^1(B(0,\sqrt{t}))}}{h_0(\sqrt{t})}+t^{\frac{N}{2p'}}\|\phi\|_{L^{p,\sigma}(B(0,\delta\sqrt{t})^c)}\right]\\
   & \le C_2t^{-\frac{N}{2}\left(1-\frac{1}{q}\right)-\frac{\alpha}{2}}
   \Gamma_{p',\sigma'}(t)\|\phi\|_{L^{p,\sigma}}
  \end{split}
  \end{equation*}
  for $\phi\in C_c({\bf R}^N)$ and $t>0$. 
  \item[{\rm (b)}] 
  Let $(p,q,\sigma,\theta)\in\Lambda$. Then there exist $C_3>0$ and $C_4>0$ such that 
  \begin{equation*}
  \begin{split}
  \|e^{-tH}\phi\|_{L^{q,\theta}}
   & \le C_3t^{-\frac{N}{2}}\frac{\|h_0\|_{L^{q,\theta}(B(0,\sqrt{t}))}}{h_0(\sqrt{t})}
  \left[\frac{\|h_0\phi\|_{L^1(B(0,\sqrt{t}))}}{h_0(\sqrt{t})}
  +t^{\frac{N}{2p'}}\|\phi\|_{L^{p,\sigma}(B(0,\sqrt{t})^c)}\right]\\
   & \le C_4t^{-\frac{N}{2}}\Gamma_{p',\sigma'}(t)\Gamma_{q,\theta}(t)\|\phi\|_{L^{p,\sigma}}
  \end{split}
  \end{equation*}
  for $\phi\in C_c({\bf R}^N)$ and $t>0$. 
\end{itemize}
\end{proposition}
\begin{proposition}
\label{Proposition:2.4}
Assume the same conditions as in Proposition~{\rm\ref{Proposition:2.3}}. 
Furthermore, assume that $h_0\in L^{p,\sigma}(B(0,1))\cap L^{p',\sigma'}(B(0,1))$ 
for some $(p,p,\sigma,\sigma)\in\Lambda$. 
Let $\|\phi\|_{p,\sigma}\le 1$ and $v_{k,i}$ be as in \eqref{eq:1.12} and set 
$$
w_{k,i}(|x|,t):=\frac{v_{k,i}(|x|,t)}{h_k(|x|)}\quad\mbox{for}\quad(x,t)\in{\bf R}^N\times(0,\infty).
$$
Then, for any $\beta\in\{0,1,2,\dots\}$, 
there exist $C>0$ and $\delta\in(0,1]$ such that 
\begin{equation}
\label{eq:2.14}
t^\beta|\partial_t^\beta w_{k,i}(x,t)|\le CM_{k,i}t^{-\frac{N}{2}}\frac{\Gamma_{p',\sigma'}(t)}{h_k(\delta\sqrt{t})}
\end{equation}
for $x\in B(0,\delta\sqrt{t})$, $t>0$, and $(k,i)\in{\mathcal K}$.
Here $M_{k,i}:=\|Q_{k,i}\|_{L^\infty({\bf S}^{N-1})}$.
\end{proposition}
Furthermore, we have:
\begin{proposition}
\label{Proposition:2.5}
Assume the same conditions as in Proposition~{\rm\ref{Proposition:2.4}}. 
Let $\|\phi\|_{p,\sigma}\le 1$ and $u_{k,i}$ be as in \eqref{eq:1.12}. 
Then, for any $\alpha\in\{0,1,\dots,m+1\}$, 
there exist $C>0$, $D>0$, and $\delta\in(0,1)$ such that
\begin{equation}
\label{eq:2.15}
|\nabla^\alpha u_{k,i}(x,t)|\le C(k+1)^Dt^{-\frac{N}{2}}\Gamma_{p',\sigma'}(t)|x|^{-\alpha}\frac{h_k(|x|)}{h_k(\delta\sqrt{t})}
\end{equation}
for $x\in B(0,\delta\sqrt{t})$, $t>0$, and $(k,i)\in{\mathcal K}$. 
\end{proposition}
{\bf Proof.}
It follows from \cite[Proposition~4.2]{IT01} that 
$$
\partial_t^\beta w_{k,i}(|x|,t)=\partial_t^\beta w_{k,i}(0,t)
+I_k[\partial_t^{\beta+1}w_{k,i}(\cdot,t)](|x|)
$$
for $(x,t)\in{\bf R}^N\times(0,\infty)$, $\beta\in\{0,1,2,\dots\}$, and $(k,i)\in{\mathcal K}$.
Repeating this relation, we have 
\begin{equation}
\label{eq:2.16}
\begin{split}
w_{k,i}(|x|,t) & =w_{k,i}(0,t)+I_k[\partial_t w_{k,i}(\cdot,t)](|x|)\\
 & =w_{k,i}(0,t)+\partial_tw_{k,i}(0,t)I_k(|x|)+I_k^2[\partial_t^2 w_{k,i}(\cdot,t)](|x|)\\
 & =\sum_{\ell=0}^{n-1}\partial_t^\ell w_{k,i}(0,t)I_k^\ell(|x|)+I_k^n[\partial_t^n w_{k,i}(\cdot,t)](|x|)
\end{split}
\end{equation}
for $(x,t)\in{\bf R}^N\times(0,\infty)$, $n\in\{1,2,\dots\}$, and $(k,i)\in{\mathcal K}$.
This implies that 
\begin{equation}
\label{eq:2.17}
u_{k,i}(x,t)=h_k(|x|)w_{k,i}(|x|,t)Q_{k,i}\left(\frac{x}{|x|}\right)
=\sum_{\ell=0}^{n-1} u_{k,i}^\ell(x,t)+R_{k,i}^n(x,t)
\end{equation}
for $(x,t)\in{\bf R}^N\times(0,\infty)$, $n\in\{1,2,\dots\}$, and $(k,i)\in{\mathcal K}$, where
\begin{equation*}
\begin{split}
 & u_{k,i}^\ell(x,t):=\partial_t^\ell w_{k,i}(0,t)J_{k,i}^\ell(|x|),\\
 & R_{k,i}^n(x,t):=h_k(|x|)I_k^n[\partial_t^n w_{k,i}(\cdot,t)](|x|)Q_{k,i}\left(\frac{x}{|x|}\right). 
\end{split}
\end{equation*}
Then, 
for any $n\in\{1,2,\dots\}$ and $\ell\in\{0,1,\dots,n-1\}$, 
by Proposition~\ref{Proposition:2.1}, Lemma~\ref{Lemma:2.1}, \eqref{eq:2.5}, and \eqref{eq:2.14} we have 
\begin{equation}
\label{eq:2.18}
\begin{split}
|\nabla^\alpha u_{k,i}^\ell(x,t)| & 
\le C(k+1)^D t^{-\frac{N}{2}-\ell}\frac{\Gamma_{p',\sigma'}(t)}{h_k(\delta\sqrt{t})}|\nabla^\alpha J_{k,i}^\ell(x)|\\
 & \le C(k+1)^D t^{-\frac{N}{2}-\ell}\frac{\Gamma_{p',\sigma'}(t)}{h_k(\delta\sqrt{t})}|x|^{2\ell-\alpha}h_k(|x|)\\
 & \le C\delta^{2\ell}(k+1)^D t^{-\frac{N}{2}}\Gamma_{p',\sigma'}(t)|x|^{-\alpha}\frac{h_k(|x|)}{h_k(\delta\sqrt{t})},\\
|\nabla^\alpha R_{k,i}^n(x,t)| & 
\le C(k+1)^D t^{-\frac{N}{2}-n}\frac{\Gamma_{p',\sigma'}(t)}{h_k(\delta\sqrt{t})}|x|^{2n-\alpha}h_k(|x|)\\
 & \le C\delta^{2n-\frac{\alpha}{2}}(k+1)^D
 t^{-\frac{N}{2}-\frac{\alpha}{2}}\Gamma_{p',\sigma'}(t)\frac{h_k(|x|)}{h_k(\delta\sqrt{t})}\quad\mbox{if}\quad \alpha\le 2n,
\end{split}
\end{equation}
for $x\in B(0,\delta\sqrt{t})$, $t>0$, and $(k,i)\in{\mathcal K}$. 
Here the constant $\delta$ is as in Proposition~\ref{Proposition:2.4} with $\beta=n$. 
Taking large enough $n\in\{0,1,2,\dots\}$ so that $2n\ge\alpha$,
by \eqref{eq:1.2} and \eqref{eq:2.17} 
we obtain inequality~\eqref{eq:2.15}
for $x\in B(0,\delta\sqrt{t})$ and $t>0$. 
Thus Proposition~\ref{Proposition:2.5} follows.
$\Box$
\section{Upper decay estimates}
In this section we study upper decay estimates of $\|\nabla^\alpha e^{-tH}\|_{(L^{p,\sigma}\to L^{q,\theta})}$ 
and prove the following theorem, which is one of the main ingredients of this paper.
\begin{theorem}
\label{Theorem:3.1}
Let $m\in\{1,2,\dots\}$ and assume conditions~{\rm ($\mbox{V}_m$)} and {\rm (N')}. 
Let $(p,q,\sigma,\theta)\in\Lambda$ and $\alpha\in\{0,1,\dots,m+1\}$.
\begin{itemize}
  \item[{\rm (a)}] There exists $C_1>0$ such that 
  $$
  \|\nabla^\alpha e^{-tH}\|_{(L^{p,\sigma}\to L^{q,\theta}(B(0,\sqrt{t})^c))}
  \le C_1t^{-\frac{N}{2}\left(1-\frac{1}{q}\right)-\frac{\alpha}{2}}
  \Gamma_{p',\sigma'}(t),\quad t>0.
  $$ 
  \item[{\rm (b)}] 
  There exists $C_2>0$ such that 
  $$
  \|\nabla^\alpha e^{-tH}\|_{(L^{p,\sigma}\to L^{q,\theta})}
  \le C_2t^{-\frac{N}{2}}\Gamma_{p',\sigma'}(t)\left[J_\alpha(t)+t^{\frac{N}{2q}-\frac{\alpha}{2}}\right],
  \quad t>0,
  $$
   where 
   $$
   J_\alpha(t):=
   \sum_{0\le k+2n\le\alpha}\sum_{i=1}^{d_k}
   t^{-n}\frac{\|\nabla^\alpha J_{k,i}^n\|_{L^{q,\theta}(B(0,\sqrt{t}))}}{h_k(\sqrt{t})}.
   $$
\end{itemize}
\end{theorem}
{\bf Proof.}
Assertion~(a) follows from Proposition~\ref{Proposition:2.3}~(a).
Then, for the proof of assertion~(b), it suffices to prove 
\begin{equation}
  \label{eq:3.1}
  \|\nabla^\alpha e^{-tH}\|_{(L^{p,\sigma}\to L^{q,\theta}(B(0,\sqrt{t})))}
  \le Ct^{-\frac{N}{2}}\Gamma_{p',\sigma'}(t)\left[J_\alpha(t)+t^{\frac{N}{2q}-\frac{\alpha}{2}}\right],
  \quad t>0.
\end{equation}

Let $\phi\in C_c({\bf R}^N)$ be such that $\|\phi\|_{p,\sigma}\le 1$ and 
set $u:=e^{-tH}\phi$. 
We use the same notations as in \eqref{eq:1.11}, \eqref{eq:1.12}, and \eqref{eq:1.13}.
For any $\ell\in\{0,1,\dots\}$, we set 
$$
[R_\ell u](x,t):=\sum_{k=\ell}^\infty\sum_{i=1}^{d_k} u_{k,i}(x,t)=u(x,t)-\sum_{k=0}^{\ell-1}\sum_{i=1}^{d_k} u_{k,i}(x,t).
$$
Let $\alpha\in\{0,1,\dots,m+1\}$ and $\delta\in(0,1)$ be small enough. 
We observe from \eqref{eq:1.2} and \eqref{eq:2.15} that
\begin{equation}
\label{eq:3.2}
\begin{split}
\left|\nabla^\alpha [R_\ell u](x,t)\right|
 & \le\sum_{k=\ell}^\infty\sum_{i=1}^{d_k}|\nabla^\alpha u_{k,i}(x,t)|\\
 & \le Ct^{-\frac{N}{2}}\Gamma_{p',\sigma'}(t)|x|^{-\alpha}
\sum_{k=\ell}^\infty(k+1)^D\frac{h_k(|x|)}{h_k(\delta\sqrt{t})}
\end{split}
\end{equation}
for $x\in B(0,\delta\sqrt{t})\setminus\{0\}$ and $t>0$. 
On the other hand, by Proposition~\ref{Proposition:2.2}~(iii), 
we find $\gamma>0$ such that 
$$
\frac{h_k(|x|)}{h_k(\delta\sqrt{t})}
=\frac{h_\ell(|x|)}{h_\ell(\delta\sqrt{t})}\frac{h_k(|x|)}{h_\ell(|x|)}\biggr/\frac{h_k(\delta\sqrt{t})}{h_\ell(\delta\sqrt{t})}
\le C\epsilon^{\left(\frac{k}{2}-\gamma\right)_+}\frac{h_\ell(|x|)}{h_\ell(\delta\sqrt{t})}
$$
for $x\in B(0,\epsilon\delta\sqrt{t})\setminus\{0\}$, $t>0$, $\epsilon\in(0,1)$, and $k\in\{\ell+1,\ell+2,\dots\}$. 
Due to the relation that $h_\ell(r/2)\asymp h_\ell(r)\asymp h_\ell(2r)$ for $r>0$ (see \eqref{eq:1.5} and \eqref{eq:1.6}),
taking small enough $\epsilon>0$ if necessary, we see that
\begin{equation}
\label{eq:3.3}
\begin{split}
\sum_{k=\ell}^\infty(k+1)^D \frac{h_k(|x|)}{h_k(\delta\sqrt{t})}
 & \le C\frac{h_\ell(|x|)}{h_\ell(\delta\sqrt{t})}
\biggr[(\ell+1)^D+\sum_{k=\ell+1}^\infty\epsilon^{\left(\frac{k}{2}-\gamma\right)_+}(k+1)^D\biggr]\\
 & \le C\frac{h_\ell(|x|)}{h_\ell(\delta\sqrt{t})}
\le C\frac{h_\ell(|x|)}{h_\ell(\sqrt{t})}
\end{split}
\end{equation}
for $x\in B(0,\epsilon\delta\sqrt{t})\setminus\{0\}$ and $t>0$. 
By \eqref{eq:3.2} and \eqref{eq:3.3} we obtain 
\begin{equation}
\label{eq:3.4}
\left|\nabla^\alpha [R_\ell u](x,t)\right|
\le Ct^{-\frac{N}{2}}\Gamma_{p',\sigma'}(t)
\frac{h_\ell^{\langle\alpha\rangle}(|x|)}{h_\ell(\sqrt{t})}
\end{equation}
for $x\in B(0,\epsilon\delta\sqrt{t})\setminus\{0\}$ and $t>0$.
Here $h_\ell^{\langle\alpha\rangle}$ is as in Subsection~2.3. 
On the other hand, 
by Proposition~\ref{Proposition:2.3}~(a) with $q=\theta=\infty$ 
we have 
\begin{equation}
\label{eq:3.5}
\begin{split}
\|\nabla^\alpha u(\cdot,t)\|_{L^{q,\theta}(B(0,\sqrt{t})\setminus B(0,\epsilon\delta\sqrt{t}))}
\le C\|\nabla^\alpha u(\cdot,t)\|_{L^\infty(B(0,\epsilon\delta\sqrt{t})^c)}t^{\frac{N}{2q}}
\le Ct^{-\frac{N}{2}+\frac{N}{2q}-\frac{\alpha}{2}}\Gamma_{p',\sigma'}(t)
\end{split}
\end{equation}
for $t>0$.

We study the behavior of $\nabla^\alpha u$ in $B(0,\epsilon\delta\sqrt{t})$ 
by using the arguments in Subsection~2.3. 
\vspace{3pt}
\newline
\underline{Step 1}: 
Consider case (${\mathcal A}_\alpha$), 
that is, $h_0(r)\not\asymp r^A$ as $r\to 0$ for $A\in\{0,1,\dots\}$ with $A\le\alpha-1$. 
Let $R_1$ be as in Subsection~2.3.
It follows from \eqref{eq:2.6} that 
\begin{equation}
\label{eq:3.6}
|\nabla^\alpha h_0(|x|)|\asymp h_0^{\langle\alpha\rangle}(|x|)>0,\quad x\in B(0,R_1)\setminus\{0\}. 
\end{equation}
Then, by  \eqref{eq:3.4} with $\ell=0$ 
we have
$$
 |\nabla^\alpha u(x,t)|=\left|\nabla^\alpha [R_0 u](x,t)\right| \le Ct^{-\frac{N}{2}}\Gamma_{p',\sigma'}(t)
\frac{h_0^{\langle\alpha\rangle}(|x|)}{h_0(\sqrt{t})}\le Ct^{-\frac{N}{2}}\Gamma_{p',\sigma'}(t)
\frac{|\nabla^\alpha h_0(|x|)|}{h_0(\sqrt{t})}
$$
for $x\in [B(0,\epsilon\delta\sqrt{t})\cap B(0,R_1)]\setminus\{0\}$ and $t>0$. 
This implies that 
\begin{equation}
\label{eq:3.7}
\|\nabla^\alpha u(\cdot,t)\|_{L^{q,\theta}(B(0,R))}
\le Ct^{-\frac{N}{2}}\frac{\Gamma_{p',\sigma'}(t)}{h_0(\sqrt{t})}\|\nabla^\alpha h_0\|_{L^{q,\theta}(B(0,R))}
\end{equation}
for $R\in(0,R_1]$ and $t>0$ if $R\le\epsilon\delta\sqrt{t}$. 

Let $R_2>1$ be such that $R_2\in(R_1,\infty)$. 
We observe from \eqref{eq:3.4} with $\ell=0$ that
$$
|\nabla^\alpha u(x,t)|=\left|\nabla^\alpha [R_0 u](x,t)\right|
\le Ct^{-\frac{N}{2}}\Gamma_{p',\sigma'}(t)|x|^{-\alpha}
\frac{h_0(|x|)}{h_0(\sqrt{t})}\le Ct^{-\frac{N}{2}}\frac{\Gamma_{p',\sigma'}(t)}{h_0(\sqrt{t})}
$$
for $x\in [B(0,\epsilon\delta\sqrt{t})\cap B(0,R_2)]\setminus B(0,R_1)$ and $t>0$. 
This together with \eqref{eq:3.6} implies that
\begin{equation}
\label{eq:3.8}
\begin{split}
\|\nabla^\alpha u(\cdot,t)\|_{L^{q,\theta}(B(0,R)\setminus B(0,R_1))}
 & \le Ct^{-\frac{N}{2}}\frac{\Gamma_{p',\sigma'}(t)}{h_0(\sqrt{t})}
\le Ct^{-\frac{N}{2}}\Gamma_{p',\sigma'}(t)\frac{\|\nabla^\alpha h_0\|_{L^{q,\theta}(B(0,R_1))}}{h_0(\sqrt{t})}\\
 & \le Ct^{-\frac{N}{2}}\Gamma_{p',\sigma'}(t)\frac{\|\nabla^\alpha h_0\|_{L^{q,\theta}(B(0,R))}}{h_0(\sqrt{t})}
\end{split}
\end{equation}
for $R\in[R_1,R_2)$ and $t>0$ if $R\le\epsilon\delta\sqrt{t}$.
Combining \eqref{eq:3.7} and \eqref{eq:3.8}, we see that
\begin{equation}
\label{eq:3.9}
\|\nabla^\alpha u(\cdot,t)\|_{L^{q,\theta}(B(0,R))}
\le Ct^{-\frac{N}{2}}\frac{\Gamma_{p',\sigma'}(t)}{h_0(\sqrt{t})}\|\nabla^\alpha h_0\|_{L^{q,\theta}(B(0,\sqrt{t}))}
\le Ct^{-\frac{N}{2}}\frac{\Gamma_{p',\sigma'}(t)}{h_0(\sqrt{t})}J_\alpha(t)
\end{equation}
for $R\in(0,R_2)$ and $t>0$ if $R<\epsilon\sqrt{t}$.
\vspace{5pt}
\newline
\underline{Step 2}: 
Consider case (${\mathcal B}_\alpha$), that is, 
$h_0(r)\asymp r^A$ as $r\to 0$, where $A\in\{0,1,\dots\}$ with $A\le\alpha-1$. 
Then $\alpha\ge 1$ and $\lambda_1=\omega_A\ge 0$ (see \eqref{eq:1.4}). Let $R_2$ as in the above.
Since $A_{1,\alpha}\ge A^+_{\omega_\alpha}=\alpha$, it follows that
$$
h_\alpha^{\langle\alpha\rangle}(|x|)\le C|x|^{A_{1,\alpha}-\alpha}\le C\min\{R,\sqrt{t}\}^{A_{1,\alpha}-\alpha},
\quad
x\in B(0,R)\cap B(0,\sqrt{t}),
$$
for $R\in(0,R_2)$.
Let $T>1$. By \eqref{eq:1.5} and \eqref{eq:1.6} we have 
\begin{equation}
\label{eq:3.10}
h_\alpha(\sqrt{t})\ge 
\left\{
\begin{array}{ll}
Ct^{\frac{A_{1,\alpha}}{2}}
 & \quad\mbox{for $0<t<T$},\vspace{5pt}\\
Ch_0(\sqrt{t}) & \quad\mbox{for $t>T$ in case $({\mathcal A}'_\alpha)$},\vspace{5pt}\\
Ct^{\frac{\alpha}{2}} & \quad\mbox{for $t>T$ in case $({\mathcal B}'_\alpha$)}.\vspace{5pt}
\end{array}
\right.
\end{equation}
These together with \eqref{eq:2.8}, \eqref{eq:2.12}, and \eqref{eq:3.4} imply that
\begin{equation}
\label{eq:3.11}
\left\|\nabla^\alpha [R_\alpha u](\cdot,t)\right\|_{L^{q,\theta}(B(0,R))}
\le Ct^{-\frac{N}{2}}\Gamma_{p',\sigma'}(t)
\frac{\|h_\alpha^{\langle\alpha\rangle}\|_{L^{q,\theta}(B(0,R))}}{h_\alpha(\sqrt{t})}
\le Ct^{-\frac{N}{2}}\Gamma_{p',\sigma'}(t)\Psi_T(t)
\end{equation}
for $R\in(0,R_2)$ and $t>0$ if $R<\epsilon\delta\sqrt{t}$. 
Here
\begin{equation}
\label{eq:3.12}
\Psi_T(t):=
 \left\{
 \begin{array}{ll}
 t^{\frac{N}{2q}-\frac{\alpha}{2}} & \mbox{if $0<t\le T$},\vspace{5pt}\\
 h_0(\sqrt{t})^{-1}& \mbox{if $t>T$ in case~$({\mathcal A}_\alpha')$},\vspace{5pt}\\
 t^{\frac{N}{2q}-\frac{\alpha}{2}} & \mbox{if $t>T$ in case~$({\mathcal B}_\alpha')$}.\vspace{5pt}
 \end{array}
 \right.
\end{equation}
Let $k\in\{0,1,2,\dots,\alpha-1\}$ and $i\in\{1,\dots,d_k\}$. 
Let $n_k\in\{0,1,2,\dots\}$ be such that 
\begin{equation}
\label{eq:3.13}
\frac{\alpha-k}{2}\le n_k<1+\frac{\alpha-k}{2}.
\end{equation} 
By \eqref{eq:2.17} and \eqref{eq:2.18} we have 
\begin{equation}
\label{eq:3.14}
\begin{split}
\left|\nabla^\alpha u_{k,i}(x,t)\right|
 & \le\sum_{\ell=0}^{n_k-1}|\nabla^\alpha u_{k,i}^\ell(x,t)|+|\nabla^\alpha R_{k,i}^{n_k}(x,t)|\\
 & \le C\sum_{\ell=0}^{n_k-1}t^{-\frac{N}{2}-\ell}\frac{\Gamma_{p',\sigma'}(t)}{h_k(\delta\sqrt{t})}|\nabla^\alpha J_{k,i}^\ell(x)|
+Ct^{-\frac{N}{2}-n_k}\frac{\Gamma_{p',\sigma'}(t)}{h_k(\delta\sqrt{t})}|x|^{2n_k-\alpha}h_k(|x|)
\end{split}
\end{equation}
for $x\in B(0,\delta\sqrt{t})$ and $t>0$. 
On the other hand, it follows from \eqref{eq:3.13} that $2n_k-\alpha+A_{1,k}\ge 2n_k-\alpha+k\ge 0$ and 
\begin{equation}
\label{eq:3.15}
|x|^{2n_k-\alpha}h_k(|x|)\le C|x|^{2n_k-\alpha+A_{1,k}}
\le C\min\{R_2,\sqrt{t}\}^{2n_k-\alpha+A_{1,k}}
\end{equation} 
for $x\in B(0,R_2)\cap B(0,\sqrt{t})$. 
Furthermore, 
similarly to \eqref{eq:3.10}, by \eqref{eq:3.13} we have
\begin{equation}
\label{eq:3.16}
\begin{split}
t^{n_k}h_k(\delta\sqrt{t}) & \ge 
\left\{
\begin{array}{ll}
Ct^{\frac{A_{1,k}}{2}+n_k}
 & \quad\mbox{if $0<t<T$},\vspace{5pt}\\
Ct^{n_k}h_0(\sqrt{t}) & \quad\mbox{if $t>T$ in case (${\mathcal A}'_\alpha$)},\vspace{5pt}\\
Ct^{\frac{\alpha-k}{2}}h_k(\sqrt{t}) & \quad\mbox{if $t>T$ in case (${\mathcal B}'_\alpha$)},\vspace{5pt}
\end{array}
\right.\\
 & \ge
\left\{
\begin{array}{ll}
Ct^{\frac{A_{1,k}}{2}+n_k}
 & \quad\mbox{if $0<t<T$},\vspace{5pt}\\
Ct^{n_k}h_0(\sqrt{t}) & \quad\mbox{if $t>T$ in case $({\mathcal A}'_\alpha$)},\vspace{5pt}\\
Ct^{\frac{\alpha}{2}}\qquad\qquad & \quad\mbox{if $t>T$ in case (${\mathcal B}'_\alpha$)},\vspace{5pt}
\end{array}
\right.
\end{split}
\end{equation}
for $t>0$.
By~\eqref{eq:3.12}, \eqref{eq:3.14}, \eqref{eq:3.15}, and \eqref{eq:3.16} 
we obtain
\begin{equation}
\label{eq:3.17}
 \|\nabla^\alpha u_{k,i}(\cdot,t)\|_{L^{q,\theta}(B(0,R))}\le Ct^{-\frac{N}{2}}\Gamma_{p',\sigma'}(t)
 \left[\sum_{\ell=0}^{n_k-1}t^{-\ell}\frac{\|\nabla^\alpha J_{k,i}^\ell\|_{L^{q,\theta}(B(0,\sqrt{t}))}}{h_k(\sqrt{t})}
+C\Psi_T(t)\right]
\end{equation} 
for $R\in(0,R_2)$ and $t>0$ if $0<R<\epsilon\delta\sqrt{t}$. 
Therefore,
thanks to \eqref{eq:3.11}, \eqref{eq:3.13}, and \eqref{eq:3.17}, 
we observe that
\begin{equation}
\label{eq:3.18}
\begin{split}
\|\nabla^\alpha u(\cdot,t)\|_{L^{q,\theta}(B(0,R))}
 & \le\sum_{k=0}^{\alpha-1}\sum_{i=1}^{d_k}\|\nabla^\alpha u_{k,i}(\cdot,t)\|_{L^{q,\theta}(B(0,R))}
 +\|\nabla^\alpha[R_\alpha u](\cdot,t)\|_{L^{q,\theta}(B(0,R))}\\
 & \le Ct^{-\frac{N}{2}}\Gamma_{p',\sigma'}(t)
 \biggr[
\sum_{0\le k+2\ell\le\alpha}\sum_{k=1}^{d_k}t^{-\ell}\frac{\|\nabla^\alpha J_{k,i}^\ell\|_{L^{q,\theta}(B(0,\sqrt{t}))}}{h_k(\sqrt{t})}
+\Psi_T(t)\biggr]\\
 & \le Ct^{-\frac{N}{2}}\Gamma_{p',\sigma'}(t)
 \biggr[J_\alpha(t)+\Psi_T(t)\biggr]
\end{split}
\end{equation}
for $R\in(0,R_2)$ and $t>0$ if $R<\epsilon\delta\sqrt{t}$. 
\vspace{3pt}

Combining \eqref{eq:3.5}, \eqref{eq:3.9}, and  \eqref{eq:3.18}, 
in both of cases~(${\mathcal A}_\alpha$) and (${\mathcal B}_\alpha$),
we see that 
\begin{equation}
\label{eq:3.19}
 \|\nabla^\alpha u(\cdot,t)\|_{L^{q,\theta}(B(0,R))}\\
 \le Ct^{-\frac{N}{2}}\Gamma_{p',\sigma'}(t)
 \biggr[J_\alpha(t)+\Psi_T(t)\biggr]
\end{equation}
for $R\in(0,R_2)$ and $t>0$ if $R<\epsilon\delta\sqrt{t}$. 
This implies that inequality~\eqref{eq:3.1} holds for $t\in(0,T]$. 
It remains to prove inequality~\eqref{eq:3.1} for $t\in(T,\infty)$.
\vspace{5pt}
\newline
\underline{Step 3}: 
Consider case (${\mathcal A}'_\alpha$), that is, 
$h_0(r)\not\asymp r^A$ as $r\to \infty$ for $A\in\{0,1,\dots\}$ with $A\le\alpha-1$. 
Taking large enough $R_2>0$ if necessary, 
we see that $R_2\ge R_3+1$, where $R_3$ is as in Subsection~2.3. 
We can assume, without loss of generality, that $\epsilon\delta\sqrt{T}>R_2$. 
By \eqref{eq:2.10} and \eqref{eq:3.4} with $\ell=0$, 
applying the same argument as in Step~1 (see also \eqref{eq:3.7}), we have 
\begin{equation}
\label{eq:3.20}
\begin{split}
\|\nabla^\alpha u(\cdot,t)\|_{L^{q,\theta}(B(0,\epsilon\delta\sqrt{t})\setminus B(0,R_2))}
 & \le Ct^{-\frac{N}{2}}\Gamma_{p',\sigma'}(t)
\frac{\|\nabla^\alpha h_0\|_{L^{q,\theta}(B(0,\sqrt{t}))}}{h_0(\sqrt{t})}\\
 & \le Ct^{-\frac{N}{2}}\Gamma_{p',\sigma'}(t)J_\alpha(t),\quad t>T.
\end{split}
\end{equation}
Since $|\nabla^\alpha h_0(x)|\asymp |x|^{-\alpha}h_0(|x|)>0$ in $B(0,R_3)^c$, 
by \eqref{eq:3.19} we have 
\begin{equation}
\label{eq:3.21}
\begin{split}
\|\nabla^\alpha u(\cdot,t)\|_{L^{q,\theta}(B(0,R_2))}
  & \le Ct^{-\frac{N}{2}}\Gamma_{p',\sigma'}(t)
 \biggr[J_\alpha(t)+\frac{\|\nabla^\alpha h_0\|_{L^{q,\theta}(B(0,R_2)\setminus B(0,R_3))}}{h_0(\sqrt{t})}\biggr]\\
  & \le Ct^{-\frac{N}{2}}\Gamma_{p',\sigma'}(t)J_\alpha(t),\quad t>T.
\end{split}
\end{equation}
By \eqref{eq:3.5}, \eqref{eq:3.20}, and \eqref{eq:3.21} we see that
$$
\|\nabla^\alpha u(\cdot,t)\|_{L^{q,\theta}(B(0,\sqrt{t}))}
\le Ct^{-\frac{N}{2}}\Gamma_{p',\sigma'}(t)
\biggr[J_\alpha(t)+t^{\frac{N}{2q}-\frac{\alpha}{2}}\biggr],\quad t>T.
$$
This implies that inequality~\eqref{eq:3.1} holds for $t>T$ in case (${\mathcal A}'_\alpha$).
\vspace{5pt}
\newline
\underline{Step 4}: 
Consider case (${\mathcal B}'_\alpha$), that is, 
$h_0(r)\asymp r^A$ as $r\to \infty$, where $A\in\{0,1,\dots\}$ with $A\le\alpha-1$.
Since 
$A_{2,\alpha}\ge A^+_{\omega_\alpha}=\alpha$, 
similarly to \eqref{eq:3.11}, by \eqref{eq:3.4} we have 
\begin{equation}
\label{eq:3.22}
\left\|\nabla^\alpha [R_\alpha u](\cdot,t)\right\|_{L^{q,\theta}(B(0,\sqrt{t})\setminus B(0,R_2))}
\le Ct^{-\frac{N}{2}}\Gamma_{p',\sigma'}(t)t^{\frac{N}{2q}-\frac{\alpha}{2}},\quad t>T. 
\end{equation} 

Let $k\in\{0,1,2,\dots,\alpha-1\}$ and $i\in\{1,\dots,d_k\}$. 
Let $n_k\in\{0,1,2,\dots\}$ be as in \eqref{eq:3.13}. 
It follows from \eqref{eq:1.6} that $2n_k-\alpha+A_{2,k}\ge 2n_k-\alpha+k\ge 0$ and 
$$
|x|^{2n_k-\alpha}h_k(|x|)\le C|x|^{2n_k-\alpha+A_{2,k}}(\log|x|)^{B_k}
\le Ct^{n_k-\frac{\alpha}{2}}h_k(\sqrt{t})
$$
for $x\in B(0,\sqrt{t})\cap B(0,R_2)^c$. 
Then, similarly to \eqref{eq:3.17}, 
we obtain 
\begin{equation}
\label{eq:3.23}
\begin{split}
 & \|\nabla^\alpha u_{k,i}(\cdot,t)\|_{L^{q,\theta}(B(0,\epsilon\delta\sqrt{t})\setminus B(0,R_2))}\\
 & \le Ct^{-\frac{N}{2}}\Gamma_{p',\sigma'}(t)
 \left[\sum_{\ell=0}^{n_k-1}t^{-\ell}\frac{\|\nabla^\alpha J_{k,i}^\ell\|_{L^{q,\theta}(B(0,\sqrt{t}))}}{h_k(\sqrt{t})}
+t^{\frac{N}{2q}-\frac{\alpha}{2}}\right],\quad t>T.
\end{split}
\end{equation}
Similarly to \eqref{eq:3.18}, 
by \eqref{eq:3.13}, \eqref{eq:3.22}, and \eqref{eq:3.23} we obtain 
\begin{equation}
\label{eq:3.24}
\begin{split}
 & \|\nabla^\alpha u(\cdot,t)\|_{L^{q,\theta}(B(0,\epsilon\delta\sqrt{t})\setminus B(0,R_2))}\\
 & \le\sum_{k=0}^{\alpha-1}\sum_{i=1}^{d_k}\|\nabla^\alpha u_{k,i}(\cdot,t)\|_{L^{q,\theta}(B(0,\epsilon\delta\sqrt{t})\setminus B(0,R_2))}
 +\|\nabla^\alpha[R_\alpha u](\cdot,t)\|_{L^{q,\theta}(B(0,\epsilon\delta\sqrt{t})\setminus B(0,R_2))}\\
 & \le Ct^{-\frac{N}{2}}\Gamma_{p',\sigma'}(t)
 \biggr[
\sum_{0\le k+2\ell\le\alpha}\sum_{i=1}^{d_k}t^{-\ell}\frac{\|\nabla^\alpha J_{k,i}^\ell\|_{L^{q,\theta}(B(0,\sqrt{t}))}}{h_k(\sqrt{t})}
+t^{\frac{N}{2q}-\frac{\alpha}{2}}\biggr]\\
 & \le Ct^{-\frac{N}{2}}\Gamma_{p',\sigma'}(t)
 \biggr[J_\alpha(t)+t^{\frac{N}{2q}-\frac{\alpha}{2}}\biggr],\quad t>T.
\end{split}
\end{equation}
Combining \eqref{eq:3.5}, \eqref{eq:3.19}, and \eqref{eq:3.24}, 
we obtain 
$$
\|\nabla^\alpha u(\cdot,t)\|_{L^{q,\theta}(B(0,\sqrt{t}))}
\le Ct^{-\frac{N}{2}}\Gamma_{p',\sigma'}(t)
 \biggr[J_\alpha(t)+t^{\frac{N}{2q}-\frac{\alpha}{2}}\biggr],\quad t>T.
$$
This implies that inequality~\eqref{eq:3.1} holds for $t>T$ in case (${\mathcal B}'_\alpha$).
Then assertion~(b) follows, and the proof of Theorem~\ref{Theorem:3.1} is complete.
$\Box$
\begin{remark}
\label{Remark:3.1} Let $H=-\Delta$. 
Then, for any $k$, $n\in\{0,1,2,\dots\}$, we find $C_{k,n}>0$ such that 
$$
I_k^n(|x|)=C_{k,n}|x|^{2n},\quad x\in{\bf R}^N.
$$
Furthermore, $h_k(|x|)=|x|^k$ and $h_k(|x|)Q_{k,i}(x/|x|)$ is a homogeneous polynomial of degree $k$. 
These mean that 
$$
|\nabla^\alpha J_{k,i}^n(x)|\le C|x|^{k+2n-\alpha}\quad\mbox{if}\quad \alpha\le k+2n,
\quad
|\nabla^\alpha J_{k,i}^n(x)|=0\quad\mbox{if}\quad \alpha>k+2n,
$$
for $x\in{\bf R}^N$. 
Then, by Theorem~{\rm\ref{Theorem:3.1}}, 
for any $(p,q,\sigma,\theta)\in\Lambda$ and $\alpha\in\{0,1,2,\dots\}$, 
we obtain the well-known decay estimate for $\nabla^\alpha e^{t\Delta}$,
$$
\|\nabla^\alpha e^{t\Delta}\|_{(L^{p,\sigma}\to L^{q,\theta})}\le Ct^{-\frac{N}{2}\left(\frac{1}{p}-\frac{1}{q}\right)-\frac{\alpha}{2}},\quad t>0.
$$
\end{remark}
%
\section{Lower decay estimates}
In this section we study lower decay estimates of $\|\nabla^\alpha e^{-tH}\|_{(L^{p,\sigma}\to L^{q,\theta})}$. 
Let 
\begin{equation*}
\begin{split}
L^{p,\sigma}_r:= & \{f\in L^{p,\sigma}\,:\,\mbox{$f$ is radially symmetric in ${\bf R}^N$}\},\\
L^{p,\sigma}_{k,i}:= & \{fQ_{k,i}\in L^{p,\sigma}\,:\,\mbox{$f$ is radially symmetric in ${\bf R}^N$}\},
\end{split}
\end{equation*}
where $(k,i)\in{\mathcal K}$ and $[fQ_{k,i}](x)=f(x)Q_{k.,i}(x/|x|)$ for $x\in{\bf R}^N\setminus\{0\}$. 
%
\begin{theorem}
\label{Theorem:4.1}
Let $m\in\{1,2,\dots\}$ and assume conditions~{\rm ($\mbox{V}_m$)} and {\rm (N')}. 
Let $(p,q,\sigma,\theta)\in\Lambda$, $\alpha\in\{0,1,\dots,m+1\}$, and $k\in\{0,1,2,\dots\}$. 
Then there exist $C>0$ and $\delta\in(0,1)$ such that 
\begin{align}
\label{eq:4.1}
 &  \|\partial_r^\alpha e^{-tH_k}\|_{(L^{p,\sigma}\to L^{q,\theta}(E))}
\ge t^{-\frac{N}{2}}\frac{\Gamma^k_{p',\sigma'}(t)}{h_k(\sqrt{t})}
 \biggr[C^{-1}\|\partial_r^\alpha h_k\|_{L^{q,\theta}(E)}-Ct^{-1}\|h_k^{\langle \alpha-2\rangle}\|_{L^{q,\theta}(E)}\biggr]_+,\\
\label{eq:4.2}
 &  \|\nabla^\alpha e^{-tH}\|_{(L^{p,\sigma}_{k,i}\to L^{q,\theta}(E))}
\ge t^{-\frac{N}{2}}\frac{\Gamma^k_{p',\sigma'}(t)}{h_k(\sqrt{t})}
 \biggr[C^{-1}\|\nabla^\alpha J_{k,i}\|_{L^{q,\theta}(E)}-Ct^{-1}\|h_k^{\langle \alpha-2\rangle}\|_{L^{q,\theta}(E)}\biggr]_+,
\end{align}
for $t>0$ and measurable sets $E\subset B(0,\delta\sqrt{t})$. 
\end{theorem}
{\bf Proof.}
Let $t>0$ and fix it.  Let $k\in\{0,1,2,\dots\}$. 
Assume that $h_k\in L^{p',\sigma'}(B(0,1))$. 
Then we find $c>0$ and a radially symmetric nonnegative function $\phi\in C_c({\bf R}^N)$ such that 
$1/2\le\|\phi\|_{L^{p,\sigma}}\le 1$ and  
\begin{equation}
\label{eq:4.3}
\int_{B(0,\sqrt{t})}h_k(y)\phi(y)\,dy
\ge c\|h_k\|_{L^{p',\sigma'}(B(0,\sqrt{t}))}.
\end{equation}
Here the constant $c$ is independent of $t$. 
Let 
\begin{equation}
\label{eq:4.4}
v(|x|,\tau):=[e^{-\tau H_k}\phi](|x|),
\qquad
w(|x|,\tau):=\frac{[e^{-\tau H_k}\phi](|x|)}{h_k(|x|)}. 
\end{equation}
Then it follows that 
\begin{equation}
\label{eq:4.5}
[e^{-\tau H}\phi_{k,i}](x)=v(|x|,\tau)Q_{k,i}\left(\frac{x}{|x|}\right),
\quad\mbox{where}\quad
\phi_{k,i}(x):=\phi(|x|)Q_{k,i}\left(\frac{x}{|x|}\right).
\end{equation}
For any $n\in\{0,1,2,\dots\}$, set 
$$
\nu_{k,n}(|{\bf x}|):=|{\bf x}|^{-n}\nu_k(|{\bf x}|)=|{\bf x}|^{-n}h_k(|{\bf x}|)^2,
\qquad
{\bf w}({\bf x},\tau):=w(|{\bf x}|,\tau),
$$
for $({\bf x},\tau)\in{\bf R}^{N+n}\times(0,\infty)$. 
Then ${\bf w}$ satisfies 
\begin{equation}
\label{eq:4.6}
\begin{split}
\partial_\tau{\bf w} =\partial_\tau w & =\frac{1}{r^{N-1}\nu_k(r)}\partial_r(r^{N-1}\nu_k(r)\partial_r w)\\
 & =\frac{1}{r^{N+n-1}r^{-n}\nu_k(r)}\partial_r(r^{N+n-1}r^{-n}\nu_k(r)\partial_r w)\\
 & =\frac{1}{\nu_{k,n}(|{\bf x}|)}\mbox{div}_{N+n}\,(\nu_{k,n}(|{\bf x}|)\nabla_{N+n}{\bf w})
\end{split}
\end{equation}
in ${\bf R}^{N+n}\times(0,\infty)$, where $r=|{\bf x}|$.

On the other hand, 
it follows from condition~(N') and \eqref{eq:1.4} that 
$2A_{1,k}\ge -N+2$ and $2A_{2,k}\ge 2A_{2,0}>-N$. 
Taking large enough $n$ if necessary, we see that
$$
-N-n<2A_{1,k}-n<N+n,
\quad
-N-n<2A_{2,k}-n<N+n.
$$
These imply that $\nu_{k,n}$ is an $A_2$-weight in ${\bf R}^{N+n}$. 
Then we apply \cite[Theorem~1.2]{IKO} (see also~\cite{CUR}) to obtain 
$$
{\bf w}(|{\bf x}|,t)\ge C\int_{{\bf R}^{N+n}}
\frac{{\bf w}(|{\bf y}|,0)}{\sqrt{\omega_{n,k}({\bf B}({\bf x},\sqrt{t}))}\sqrt{\omega_{n,k}({\bf B}({\bf y},\sqrt{t}))}}
\exp\left(-\frac{|{\bf x}-{\bf y}|^2}{Ct}\right)\nu_{k,n}(|{\bf y}|)\,d{\bf y}
$$
for $({\bf x},t)\in{\bf R}^{N+n}\times(0,\infty)$. 
Here 
$$
{\bf B}({\bf x},r):=\{{\bf y}\in{\bf R}^{N+n}\,:\,|{\bf x}-{\bf y}|_{N+n}<r\},
\quad
\omega_{n,k}({\bf B}({\bf x},r)):=\int_{{\bf B}({\bf x},r)}\nu_{k,n}({\bf y})\,d{\bf y}.
$$
It follows from Proposition~\ref{Proposition:2.2}~(b) that
$$
\omega_{n,k}({\bf B}(0,r))\le C\int_0^r \nu_k(r)r^{N-1}\,dr\le C\int_{B(0,r)}\nu_k(|x|)\,dx\le Cr^Nh_k^2(r),\quad r>0.
$$
Then, recalling ${\bf w}({\bf x},0)=\phi(r)/h_k(r)$ with $r=|{\bf x}|$, 
we have 
\begin{equation*}
\begin{split}
w(x,t)={\bf w}({\bf x},t)
 & \ge C\int_{{\bf B}(0,\sqrt{t})}
\frac{{\bf w}(|{\bf y}|,0)}{\sqrt{\omega_{n,k}({\bf B}({\bf x},\sqrt{t}))}
\sqrt{\omega_{n,k}({\bf B}({\bf y},\sqrt{t}))}}\nu_{k,n}({\bf y})\,d{\bf y}\\
 & \ge C\int_{{\bf B}(0,\sqrt{t})}
 \frac{{\bf w}({\bf y},0)\nu_{k,n}({\bf y})}{\sqrt{\omega_{n,k}({\bf B}(0,2\sqrt{t}))}\sqrt{\omega_{n,k}({\bf B}(0,2\sqrt{t}))}}\,d{\bf y}\\
 & \ge Ct^{-\frac{N}{2}}\nu_k(\sqrt{t})^{-1}
 \int_0^{\sqrt{t}}\frac{\phi(r)}{h_k(r)}r^{-n}\nu_k(r)r^{N+n-1}\,dr\\
 & \ge Ct^{-\frac{N}{2}}\nu_k(\sqrt{t})^{-1} \int_0^{\sqrt{t}}\phi(r)h_k(r)r^{N-1}\,dr
\end{split}
\end{equation*}
for $x\in B(0,\sqrt{t})$. 
This together with \eqref{eq:4.3} implies that
\begin{equation}
\label{eq:4.7}
w(x,t)
\ge C^{-1}t^{-\frac{N}{2}}h_k(\sqrt{t})^{-2}\int_{B(0,\sqrt{t})}h_k(y)\phi(y)\,dy
\ge C^{-1}t^{-\frac{N}{2}}\frac{\Gamma^k_{p',\sigma'}(t)}{h_k(\sqrt{t})}
\end{equation}
for $x\in B(0,\sqrt{t})$. 

On the other hand, 
similarly to \eqref{eq:2.16}, 
taking small enough $\delta\in(0,1)$, 
we obtain 
\begin{equation}
\label{eq:4.8}
w(|x|,t)=\sum_{\ell=0}^{n-1}\partial_t^\ell w(0,t)I_k^\ell(|x|)+I_k^n[\partial_t^n w(\cdot,t)](|x|),
\quad x\in{\bf R}^N,
\end{equation}
where $n=1,2,\dots$. 
Applying Proposition~\ref{Proposition:2.4} with $H$ replaced by $H_k$, 
for any $\beta\in\{0,1,2,\dots\}$, 
we have 
\begin{equation}
\label{eq:4.9}
t^\beta|\partial_t^\beta w(x,t)|\le Ct^{-\frac{N}{2}}\frac{\Gamma^k_{p',\sigma'}(t)}{h_k(\delta\sqrt{t})}
\le Ct^{-\frac{N}{2}}\frac{\Gamma^k_{p',\sigma'}(t)}{h_k(\sqrt{t})}
\end{equation}
for $x\in B(0,\delta\sqrt{t})$. 
Then, taking small enough $\epsilon\in(0,1)$ and applying Lemma~\ref{Lemma:2.1}, 
by \eqref{eq:4.7}, \eqref{eq:4.8} with $n=1$ and \eqref{eq:4.9} 
we see that 
\begin{equation}
\label{eq:4.10}
\begin{split}
w(0,t) & \ge w(|x|,t)-I_k[\partial_tw(\cdot,t)](|x|)\\
 & \ge C^{-1}t^{-\frac{N}{2}}\frac{\Gamma^k_{p',\sigma'}(t)}{h_k(\sqrt{t})}
-Ct^{-\frac{N}{2}-1}\frac{\Gamma^k_{p',\sigma'}(t)}{h_k(\sqrt{t})}|x|^2
\ge C^{-1}t^{-\frac{N}{2}}\frac{\Gamma^k_{p',\sigma'}(t)}{h_k(\sqrt{t})}
\end{split}
\end{equation}
for $x\in B(0,\epsilon\delta\sqrt{t})$. 
By \eqref{eq:4.8}, \eqref{eq:4.9}, and \eqref{eq:4.10}, 
taking large enough $n$ if necessary, 
we apply Proposition~\ref{Proposition:2.1} and Lemma~\ref{Lemma:2.1} 
to obtain   
\begin{equation}
\label{eq:4.11}
\begin{split}
|\partial_r^\alpha v(|x|,t)| & =|\partial_r^\alpha [h_k(|x|)w(|x|,t)|\\
 & \ge |w(0,t)||\partial_r^\alpha h_k(|x|)|-C\sum_{\ell=1}^nt^{-\frac{N}{2}-\ell}\frac{\Gamma^k_{p',\sigma'}(t)}{h_k(\sqrt{t})}|x|^{2\ell-\alpha}h_k(|x|)\\
 & \ge t^{-\frac{N}{2}}\frac{\Gamma^k_{p',\sigma'}(t)}{h_k(\sqrt{t})}
 \biggr[C^{-1}|\partial_r^\alpha h_k(|x|)|-Ct^{-1}|x|^{2-\alpha}h_k(|x|)\biggr]
\end{split}
\end{equation}
for $x\in B(0,\delta\sqrt{t})\setminus\{0\}$. 
This implies that 
\begin{equation}
\label{eq:4.12}
\begin{split}
 & \|\partial_r^\alpha e^{-tH_k}\|_{(L^{p,\sigma}\to L^{q,\theta}(E))}
 \ge\frac{\|\partial_r^\alpha v(\cdot,t)\|_{L^{q,\theta}(E)}}{\|\phi\|_{L^{p,\sigma}}}\\
 & \qquad\quad
\ge t^{-\frac{N}{2}}\frac{\Gamma^k_{p',\sigma'}(t)}{h_k(\sqrt{t})}
\biggr[C^{-1}\|\partial_r^\alpha h_k\|_{L^{q,\theta}(E)}-Ct^{-1}\|h_k^{\langle \alpha-2\rangle}\|_{L^{q,\theta}(E)}\biggr]_+
\end{split}
\end{equation}
for measurable sets $E\subset B(0,\delta\sqrt{t})$. 
Thus inequality~\eqref{eq:4.1} holds when $h_k\in L^{p',\sigma'}(B(0,1))$. 

On the other hand, 
it follows from \eqref{eq:4.4}, \eqref{eq:4.5}, and \eqref{eq:4.8} that
$$
[e^{-tH}\phi_{k,i}](x)
=\sum_{\ell=0}^{n-1}\partial_t^\ell w(0,t)J_{k,i}^\ell(|x|)+Q_{k,i}(x/|x|)h_k(|x|)I_k^n[\partial_t^n w(\cdot,t)](|x|),
\quad x\in{\bf R}^N\setminus\{0\}.
$$
Similarly to \eqref{eq:4.11}, 
we have 
\begin{equation*}
\begin{split}
|\nabla^\alpha[e^{-tH}\phi_{k,i}](x)|
 & \ge|w(0,t)||\nabla^\alpha J_{k,i}^0(x)|
-C\sum_{\ell=1}^nt^{-\frac{N}{2}-\ell}\frac{\Gamma^k_{p',\sigma'}(t)}{h_k(\sqrt{t})}|x|^{2\ell-\alpha}h_k(|x|)\\
 & \ge t^{-\frac{N}{2}}\frac{\Gamma^k_{p',\sigma'}(t)}{h_k(\sqrt{t})}
 \biggr[C^{-1}|\nabla^\alpha J_{k,i}^0(x)|-Ct^{-1}|x|^{2-\alpha}h_k(|x|)\biggr]
\end{split}
\end{equation*}
for $x\in B(0,\delta\sqrt{t})\setminus\{0\}$. 
Then, similarly to \eqref{eq:4.12}, we obtain inequality~\eqref{eq:4.2} when $h_k\in L^{p',\sigma'}(B(0,1))$. 

If $h_k\not\in L^{p',\sigma'}(B(0,1))$, 
then we approximate the potential~$V$ by bounded radially symmetric potentials 
and apply the above arguments. 
Then we deduce that inequalities~\eqref{eq:4.1} and \eqref{eq:4.2} hold 
with $\Gamma_{p',\sigma'}^k(t)=\infty$ when $h_k\not\in L^{p',\sigma'}(B(0,1))$. 
Thus Theorem~\ref{Theorem:4.1} follows.
$\Box$\vspace{5pt}

We apply Theorem~\ref{Theorem:4.1} to obtain the following theorem.
\begin{theorem}
\label{Theorem:4.2}
Let $m\in\{1,2,\dots\}$ and assume conditions~{\rm ($\mbox{V}_m$)} and {\rm (N')}. 
\begin{itemize}
  \item[{\rm (a)}] 
   Let $\alpha\in\{0,1,\dots,m+1\}$ and $k\in\{0,1,2,\dots\}$. Assume that 
  \begin{equation}
  \label{eq:4.13}
  C^{-1}r^{-\alpha}h_k(r)\le|\partial_r^\alpha h_k(r)|\le Cr^{-\alpha}h_k(r)\quad\mbox{for}\quad R_1<r<R_2,
  \end{equation}
  for some $0\le R_1<R_2\le\infty$. 
  Then there exist $C_1>0$ and $\delta_1\in(0,1)$ such that 
  $$
  \|\partial_r^\alpha e^{-tH_k}\|_{(L^{p,\sigma}_r\to L^{q,\theta}(E))}
  \ge C_1^{-1}t^{-\frac{N}{2}}\Gamma^k_{p',\sigma'}(t)
  \frac{\|\partial_r^\alpha h_k\|_{L^{q,\theta}(E)}}{h_k(\sqrt{t})}
  $$
  for measurable sets $E\subset\{x\in B(0,\delta_1\sqrt{t})\,:\,R_1<|x|<R_2\}$ and $t>0$. 
  \item[{\rm (b)}] 
  Let $\alpha\in\{0,1,\dots,m+1\}$. Then there exist $C_2>0$ and $\delta_2$, $\delta_3\in(0,1)$ with $\delta_2<\delta_3$ such that
  \begin{equation}
  \label{eq:4.14}
  \|\partial_r^\alpha e^{-tH}\|_{(L^{p,\sigma}\to L^{q,\theta}(E))}
  \ge Ct^{-\frac{N}{2p}-\frac{\alpha}{2}}|E|^{\frac{N}{q}}
  \end{equation}
  for measurable sets $E\subset B(0,\delta_3\sqrt{t})\cap B(0,\delta_2\sqrt{t})^c$ and $t>0$. 
  In particular, 
  \begin{equation}
  \label{eq:4.15}
  \|\nabla^\alpha e^{-tH}\|_{(L^{p,\sigma}\to L^{q,\theta})}
  \ge\|\partial_r^\alpha e^{-tH}\|_{(L^{p,\sigma}\to L^{q,\theta})}
  \ge Ct^{-\frac{N}{2}\left(\frac{1}{p}-\frac{1}{q}\right)-\frac{\alpha}{2}},\quad t>0.
  \end{equation}
  \end{itemize}
\end{theorem}
{\bf Proof.}
We prove assertion~(a). Assume \eqref{eq:4.13}. 
Let $\delta\in(0,1)$ be as in Theorem~\ref{Theorem:4.1} and $0<\epsilon\le 1$. 
It follows that 
\begin{equation}
\label{eq:4.16}
t^{-1}h_k^{\langle \alpha-2\rangle}(|x|)=t^{-1}r^{-\alpha+2}h_k(|x|)
\le(\epsilon\delta)^2|x|^{-\alpha}h_k(|x|)
\le C(\epsilon\delta)^2|\partial_r^\alpha h_k(|x|)|
\end{equation}
for measurable sets $E\subset\{x\in B(0,\epsilon\delta\sqrt{t})\,:\,R_1<|x|<R_2\}$. 
Taking small enough $\epsilon\in(0,1)$ if necessary, 
by Theorem~\ref{Theorem:4.1} and \eqref{eq:4.16} 
we have 
\begin{equation*}
\begin{split}
\|\partial_r^\alpha e^{-tH_k}\|_{(L^{p,\sigma}\to L^{q,\theta}(E))}
 & \ge t^{-\frac{N}{2}}\frac{\Gamma^k_{p',\sigma'}(t)}{h_k(\sqrt{t})}
 \biggr[C^{-1}\|\partial_r^\alpha h_k\|_{L^{q,\theta}(E)}-C\epsilon^2\delta^2\|\partial_r^\alpha h_k\|_{L^{q,\theta}(E)}\biggr]_+\\
 & \ge C^{-1}t^{-\frac{N}{2}}\frac{\Gamma^k_{p',\sigma'}(t)}{h_k(\sqrt{t})}\|\partial_r^\alpha h_k\|_{L^{q,\theta}(E)}
\end{split}
\end{equation*}
for measurable sets $E\subset\{x\in B(0,\epsilon\delta\sqrt{t})\,:\,R_1<|x|<R_2\}$ and $t>0$. 
Thus assertion~(a) follows. 

We prove assertion~(b). 
Let $\delta'\in(0,\delta)$. 
Taking large enough $k\in\{0,1,2,\dots\}$, 
by Propositions~\ref{Proposition:2.1} and \ref{Proposition:2.2} 
we see that 
\begin{equation*}
\begin{split}
 & C^{-1}r^{-\alpha}h_k(r)\le \left|\frac{d^\alpha}{dr^\alpha} h_k(r)\right|\le Cr^{-\alpha}h_k(r),
\quad\mbox{for}\quad r>0,\\
 & \Gamma_{p',\sigma'}^k(t)\ge C^{-1}t^{-\frac{N}{2p'}},
\qquad
\frac{\|\partial_r^\alpha h_k\|_{L^{q,\theta}(E)}}{h_k(\sqrt{t})}
\ge C^{-1}t^{-\frac{\alpha}{2}}|E|^{\frac{1}{q}}\quad\mbox{for}\quad t>0,
\end{split}
\end{equation*}
where $E$ is a measurable set in $B(0,\epsilon\delta\sqrt{t})\cap B(0,\epsilon\delta'\sqrt{t})^c$. 
Then assertion~(a) implies that 
\begin{equation*}
\begin{split}
\|\partial_r^\alpha e^{-tH}\|_{(L^{p,\sigma}\to L^{q,\theta}(E))}
 & \ge\|\partial_r^\alpha e^{-tH}\|_{(L_{k,i}^{p,\sigma}\to L^{q,\theta}(E))}
\ge C^{-1}\|\partial_r^\alpha e^{-tH_k}\|_{(L^{p,\sigma}\to L^{q,\theta}(E))}\\
 & \ge C^{-1}t^{-\frac{N}{2p'}-\frac{\alpha}{2}}|E|^{\frac{1}{q}}
\end{split}
\end{equation*}
for measurable sets $E\subset B(0,\epsilon\delta\sqrt{t})\cap B(0,\epsilon\delta'\sqrt{t})^c$ and $t>0$, where $i\in\{1,\dots,d_k\}$.
Then we have inequality~\eqref{eq:4.14}. 
Inequality~\eqref{eq:4.15} easily follows from \eqref{eq:4.14}. 
Thus we obtain assertion~(b), and the proof of Theorem~\ref{Theorem:4.2} is complete.
$\Box$
\section{Proof of Theorem~\ref{Theorem:1.1}}
In this section, 
combining the results in the previous sections, we prove Theorem~\ref{Theorem:1.1}. 
\vspace{3pt}
\newline
{\bf Proof of Theorem~\ref{Theorem:1.1}.}
Let us consider the case $\alpha=0$. 
It follows from Theorem~\ref{Theorem:3.1} and Proposition~\ref{Proposition:2.2}~(a) 
that 
$$
\|e^{-tH}\|_{(L^{p,\sigma}\to L^{q,\theta})}
\le Ct^{-\frac{N}{2}}\Gamma_{p',\sigma'}(t)
\left[\frac{\|h_0\|_{L^{q,\theta}(B(0,\sqrt{t}))}}{h_0(\sqrt{t})}+t^{\frac{N}{2q}}\right]
\le Ct^{-\frac{N}{2}}\Gamma_{p',\sigma'}(t)\Gamma_{q,\theta}(t),\quad t>0.
$$
Let $\delta\in(0,1)$ be small enough. 
By Theorem~\ref{Theorem:4.2} we see that
\begin{equation*}
\begin{split}
\|e^{-tH}\|_{(L^{p,\sigma}\to L^{q,\theta})}
 & \ge\|e^{-tH}\|_{(L_r^{p,\sigma}\to L^{q,\theta}(B(0,\delta\sqrt{t})))}\\
 & \ge Ct^{-\frac{N}{2}}\Gamma_{p',\sigma'}(t)
\frac{\|h_0\|_{L^{q,\theta}(B(0,\delta\sqrt{t}))}}{h_0(\sqrt{t})}
\ge Ct^{-\frac{N}{2}}\Gamma_{p',\sigma'}(t)\Gamma_{q,\theta}(t),
\quad t>0.
\end{split}
\end{equation*}
Then Theorem~\ref{Theorem:1.1} follows in the case $\alpha=0$. 
So it suffices to prove Theorem~\ref{Theorem:1.1} in the case $\alpha\in\{1,2\}$. 
\vspace{3pt}
\newline
\underline{Step 1}: 
We prove that 
\begin{equation}
\label{eq:5.1}
\|\nabla^\alpha e^{-tH}\|_{(L^{p,\sigma}\to L^{q,\theta})}\le C\Phi_\alpha(t),\qquad t>0,
\end{equation}
where $\alpha\in\{1,2\}$. 
We can assume, without loss of generality, that $\Phi_\alpha(t_*)<\infty$ for some $t_*>0$. 
Then
\begin{equation}
\label{eq:5.2}
h_0\in L^{p',\sigma'}(B(0,\sqrt{t_*})),\qquad \nabla^\alpha h_0\in L^{q,\theta}(B(0,\sqrt{t_*})).
\end{equation}
These imply that 
\begin{equation}
\label{eq:5.3}
\begin{split}
 & h_0\in L^{p',\sigma'}(B(0,R)),\quad
 \nabla^\alpha h_0\in L^{q,\theta}(B(0,R)),
\quad
\nabla^\alpha J_{1,i}\in L^{q,\theta}(B(0,R)),\\
 & \Phi_\alpha(t)\le
J_\alpha(t)<\infty,
\end{split}
\end{equation}
for $R>0$ and $t>0$, where $i=1,\dots,N$. 
It follows from Theorem~\ref{Theorem:3.1} that 
\begin{equation}
\label{eq:5.4}
\begin{split}
 & \|\nabla e^{-tH}\|_{(L^{p,\sigma}\to L^{q,\theta})}\\
 &  \le Ct^{-\frac{N}{2}}\Gamma_{p',\sigma'}(t)
\biggr[\Phi_1(t)+\sum_{i=1}^N\frac{\|\nabla J_{1,i}\|_{L^{q,\theta}(B(0,\sqrt{t}))}}{h_1(\sqrt{t})}\biggr],\\
 & \|\nabla^2 e^{-tH}\|_{(L^{p,\sigma}\to L^{q,\theta})}\\
 & \le Ct^{-\frac{N}{2}}\Gamma_{p',\sigma'}(t)
 \biggr[\Phi_2(t)+t^{-1}\frac{\|\nabla^2 J_{0,1}^1\|_{L^{q,\theta}(B(0,\sqrt{t}))}}{h_0(\sqrt{t})}
 +\frac{\|\nabla^2 h_2\|_{L^{q,\theta}(B(0,\sqrt{t}))}}{h_2(\sqrt{t})}\biggr],
 \quad t>0.
\end{split}
\end{equation}

Following the arguments in Subsection~2.2, 
we divide the behavior of $h_0$ near $0$ into cases (${\mathcal A}_2$) and (${\mathcal B}_2$). 
In case~(${\mathcal A}_2$),  
by \eqref{eq:2.7} and \eqref{eq:5.4} we find $R_1>0$ such that 
\eqref{eq:5.1} holds for $t\in(0,\sqrt{R_1})$. 
In case~(${\mathcal B}_2$), 
for any $R_2>0$, 
by \eqref{eq:2.9} and \eqref{eq:5.4} we have \eqref{eq:5.1} for $t\in(0,\sqrt{R_2})$.
These imply that \eqref{eq:5.1} holds in $(0,T_1]$ for some $T_1>0$.  

Let $T_2\in(T_1,\infty)$ be large enough. 
By \eqref{eq:5.3} we see that 
$$
J_\alpha(t)+t^{\frac{N}{2q}-\frac{\alpha}{2}}\le C,
\quad
\Phi_\alpha(t)\ge t^{\frac{N}{2q}-\frac{\alpha}{2}}\ge C,
\quad\mbox{for}\quad t\in(T_1,T_2].
$$
Then Theorem~\ref{Theorem:3.1} implies that \eqref{eq:5.1} holds for $T_1<t\le T_2$. 

Next we we divide the behavior of $h_0$ at the space infinity 
into cases (${\mathcal A}_2'$) and (${\mathcal B}_2'$). 
Consider case~(${\mathcal A}_2'$). 
Let $R_3>0$ be as in \eqref{eq:2.10}. 
Taking large enough $T_2>0$ if necessary, 
by \eqref{eq:2.11} and \eqref{eq:5.3} we have
\begin{equation*}
\begin{split}
\sum_{i=1}^N\frac{\|\nabla J_{1,i}\|_{L^{q,\theta}(B(0,\sqrt{t}))}}{h_1(\sqrt{t})}
 & \le\sum_{i=1}^N\frac{\|\nabla J_{1,i}\|_{L^{q,\theta}(B(0,\sqrt{t})\setminus B(0,R_3))}+C}{h_1(\sqrt{t})}\\
 & \le C\frac{\|\nabla h_0\|_{L^{q,\theta}(B(0,\sqrt{t})\setminus B(0,R_3))}+1}{h_0(\sqrt{t})}
 \le C\frac{\|\nabla h_0\|_{L^{q,\theta}(B(0,\sqrt{t}))}}{h_0(\sqrt{t})},\\
t^{-1}\frac{\|\nabla^2 J_{0,1}^1\|_{L^{q,\theta}(B(0,\sqrt{t}))}}{h_0(\sqrt{t})}
 & \le t^{-1}\frac{\|\nabla^2 J_{0,1}^1\|_{L^{q,\theta}(B(0,\sqrt{t})\setminus B(0,R_3))}+C}{h_0(\sqrt{t})}\\
 & \le C\frac{\|\nabla^2 h_0\|_{L^{q,\theta}(B(0,\sqrt{t})\setminus B(0,R_3))}+1}{h_0(\sqrt{t})}
\le C\frac{\|\nabla^2 h_0\|_{L^{q,\theta}(B(0,\sqrt{t}))}}{h_0(\sqrt{t})},\\
\frac{\|\nabla^2h_2\|_{L^{q,\theta}(B(0,\sqrt{t}))}}{h_2(\sqrt{t})}
 & \le \frac{\|\nabla^2 h_2\|_{L^{q,\theta}(B(0,\sqrt{t})\setminus B(0,R_3))}+C}{h_2(\sqrt{t})}\\
 & \le C\frac{\|\nabla^2 h_0\|_{L^{q,\theta}(B(0,\sqrt{t})\setminus B(0,R_3))}+1}{h_0(\sqrt{t})}
\le C\frac{\|\nabla^2 h_0\|_{L^{q,\theta}(B(0,\sqrt{t}))}}{h_0(\sqrt{t})},
\end{split}
\end{equation*}
for large enough $t\ge T_2$. 
These together with \eqref{eq:5.4} imply that 
inequality~\eqref{eq:5.1} holds for $t\ge T_2$ in case (${\mathcal A}_2'$). 

Consider case~(${\mathcal B}_2'$). Let $R_4>0$ be large enough. 
By \eqref{eq:2.12} and \eqref{eq:2.13}, 
taking large enough $T_2$ if necessary, we have
\begin{equation}
\label{eq:A1}
\begin{split}
 & \sum_{i=1}^N\frac{\|\nabla J_{1,i}\|_{L^{q,\theta}(B(0,\sqrt{t}))}}{h_1(\sqrt{t})}
\le \sum_{i=1}^N\frac{\|\nabla J_{1,i}\|_{L^{q,\theta}(B(0,\sqrt{t})\setminus B(0,R_4))}+C}{h_1(\sqrt{t})}
\le Ct^{\frac{N}{2q}-\frac{1}{2}},\\
 & t^{-1}\frac{\|\nabla^2 J_{0,1}^1\|_{L^{q,\theta}(B(0,\sqrt{t}))}}{h_0(\sqrt{t})}
\le t^{-1}\frac{\|\nabla^2 J_{0,1}^1\|_{L^{q,\theta}(B(0,\sqrt{t})\setminus B(0,R_4))}+C}{h_0(\sqrt{t})}
\le Ct^{\frac{N}{2q}-1},\\
 & \frac{\|\nabla^2 h_2\|_{L^{q,\theta}(B(0,\sqrt{t}))}}{h_2(\sqrt{t})}
\le\frac{\|\nabla^2 h_2\|_{L^{q,\theta}(B(0,\sqrt{t})\setminus B(0,R_4))}+C}{h_2(\sqrt{t})}
\le Ct^{\frac{N}{2q}-1}+Ct^{-1},
\end{split}
\end{equation}
for $t\ge T_2$. 
By \eqref{eq:5.4} and \eqref{eq:A1} we see that inequality~\eqref{eq:5.1} holds for $t\ge T_2$ in case~(${\mathcal B}_2'$). 
Therefore we obtain \eqref{eq:5.1} for $t\ge T_2$, 
and the proof of \eqref{eq:5.1} is complete.
\vspace{3pt}
\newline
\underline{Step 2}: 
We prove 
\begin{equation}
\label{eq:5.5}
\|\nabla^\alpha e^{-tH}\|_{(L^{p,\sigma}\to L^{q,\theta})}\ge C^{-1}\Phi_\alpha(t),\qquad t>0,
\end{equation}
where $\alpha\in\{1,2\}$. 
We can assume, without loss of generality, that 
\begin{equation}
\label{eq:5.6}
\|\nabla^\alpha e^{-t_*H}\|_{(L^{p,\sigma}\to L^{q,\theta})}<\infty\quad\mbox{for some $t_*>0$}.
\end{equation}
We also prove that \eqref{eq:5.3} holds for $R>0$ and $t>0$ 
under assumption~\eqref{eq:5.6}. 
Similarly to Step~1,
we divide the behavior of $h_0$ near $0$ into cases (${\mathcal A}_2$) and (${\mathcal B}_2$). 
Let $R_i$ $(i=1,2,3,4)$ be as in the above.

Consider case~(${\mathcal A}_2$).  
By Theorem~\ref{Theorem:4.2}~(a), \eqref{eq:2.6}, and \eqref{eq:2.7} we see that 
\begin{equation}
\label{eq:5.7}
\begin{split}
 & \|\nabla^\alpha e^{-tH}\|_{(L^{p,\sigma}\to L^{q,\theta})}
\ge\|\partial^\alpha_r e^{-tH}\|_{(L^{p,\sigma}_r\to L^{q,\theta}(B(0,R)))}\\
 & \ge C^{-1}t^{-\frac{N}{2}}\Gamma_{p',\sigma'}(t)\frac{\|\partial_r^\alpha h_0\|_{L^{q,\theta}(B(0,R))}}{h_0(\sqrt{t})}\\
 & \ge C^{-1}t^{-\frac{N}{2}}\Gamma_{p',\sigma'}(t)\frac{\|\nabla^\alpha h_0\|_{L^{q,\theta}(B(0,R))}}{h_0(\sqrt{t})}
 \ge C^{-1}t^{-\frac{N}{2}}\Gamma_{p',\sigma'}(t)\frac{\|\nabla^\alpha J_{1,i}\|_{L^{q,\theta}(B(0,R))}}{h_1(\sqrt{t})}
\end{split}
\end{equation}
for $0<R<\delta\sqrt{t}$ with $R<R_1$ and $t>0$. This together with \eqref{eq:5.6} implies \eqref{eq:5.2} and \eqref{eq:5.3}. 
On the other hand, 
it follows from Proposition~\ref{Proposition:2.2}~(a) and \eqref{eq:2.6} that
\begin{equation}
\label{eq:5.8}
t^{\frac{N}{2q}-\frac{\alpha}{2}}\le Ct^{-\frac{\alpha}{2}}\frac{\|h_0\|_{L^{q,\theta}(B(0,\delta\sqrt{t}))}}{h_0(\sqrt{t})}
\le C\frac{\|h_0^{\langle\alpha\rangle}\|_{L^{q,\theta}(B(0,\delta\sqrt{t}))}}{h_0(\sqrt{t})}
\le C\frac{\|\nabla^\alpha h_0\|_{L^{q,\theta}(B(0,\delta\sqrt{t}))}}{h_0(\sqrt{t})}
\end{equation}
for $0<t\le R_1^2$. 
By \eqref{eq:5.7} and \eqref{eq:5.8} 
we see that \eqref{eq:5.5} holds for $t\in(0,R_1^2)$. 

Consider case~(${\mathcal B}_2$). 
By Theorem~\ref{Theorem:4.1} and \eqref{eq:2.8} we have
\begin{equation}
\label{eq:5.9}
\begin{split}
 & \|\nabla^\alpha e^{-tH}\|_{(L^{p,\sigma}\to L^{q,\theta})}
 \ge \|\nabla^\alpha e^{-tH}\|_{(L^{p,\sigma}_{k,i}\to L^{q,\theta}(B(0,R)))}\\
 & \ge t^{-\frac{N}{2}}\frac{\Gamma^k_{p',\sigma'}(t)}{h_k(\sqrt{t})}
\biggr[C^{-1}\|\nabla^\alpha J_{k,i}\|_{L^{q,\theta}(B(0,R))}-Ct^{-1}\|h_k^{\langle \alpha-2\rangle}\|_{L^{q,\theta}(B(0,R))}\biggr]_+\\
& \ge C^{-1}t^{-\frac{N}{2p}}\frac{\|\nabla^\alpha J_{k,i}\|_{L^{q,\theta}(B(0,R))}}{h_k(\sqrt{t})}
-Ct^{-\frac{N}{2}\left(\frac{1}{p}-\frac{1}{q}\right)-\frac{\alpha}{2}}\end{split}
\end{equation}
for $0<R<\delta\sqrt{t}$ with $R<R_2$ and $t>0$,  
where $k\in\{0,1\}$ and $i\in\{1,\dots,d_k\}$. 
This together with \eqref{eq:5.6} implies \eqref{eq:5.2} and \eqref{eq:5.3}. 

Let $\epsilon>0$ be small enough. 
By Theorem~\ref{Theorem:4.2}~(b), \eqref{eq:2.8}, and \eqref{eq:5.9} 
we obtain 
\begin{equation*}
\begin{split}
 & \|\nabla e^{-tH}\|_{(L^{p,\sigma}\to L^{q,\theta})}
 \ge (1-\epsilon)\|\nabla e^{-tH}\|_{(L^{p,\sigma}\to L^{q,\theta}(B(0,\delta\sqrt{t})))}
 +\epsilon\|\nabla e^{-tH}\|_{(L^{p,\sigma}_{0,1}\to L^{q,\theta}(B(0,\delta\sqrt{t})))}\\
 & \ge C^{-1}t^{-\frac{N}{2}\left(\frac{1}{p}-\frac{1}{q}\right)-\frac{1}{2}}
 +C^{-1}\epsilon t^{-\frac{N}{2p}}\frac{\|\nabla J_{0,1}\|_{L^{q,\theta}(B(0,\sqrt{t}))}}{h_0(\sqrt{t})}
 -C\epsilon t^{-\frac{N}{2}\left(\frac{1}{p}-\frac{1}{q}\right)-\frac{1}{2}}\\
 & \ge C^{-1}t^{-\frac{N}{2}\left(\frac{1}{p}-\frac{1}{q}\right)-\frac{1}{2}}
 +C^{-1}\epsilon t^{-\frac{N}{2p}}\frac{\|\nabla h_0\|_{L^{q,\theta}(B(0,\sqrt{t}))}}{h_0(\sqrt{t})}\\
  & \ge C^{-1}\epsilon t^{-\frac{N}{2p}}\biggr[\frac{\|\nabla h_0\|_{L^{q,\theta}(B(0,\sqrt{t}))}}{h_0(\sqrt{t})}+t^{\frac{N}{2q}-\frac{1}{2}}\biggr]
  \ge C^{-1}\epsilon\Phi_1(t),
  \quad t\in(0,R_2^2).
\end{split}
\end{equation*}
Similarly, we have
\begin{equation*}
\begin{split}
 & \|\nabla^2 e^{-tH}\|_{(L^{p,\sigma}\to L^{q,\theta})}\\
 & \ge (1-(N+1)\epsilon)\|\nabla^2 e^{-tH}\|_{(L^{p,\sigma}\to L^{q,\theta})}
 +\epsilon\|\nabla^2 e^{-tH}\|_{(L^{p,\sigma}_{0,1}\to L^{q,\theta})}
 +\epsilon\sum_{i=1}^N\|\nabla^2 e^{-tH}\|_{(L^{p,\sigma}_{1,i}\to L^{q,\theta})}\\
 & \ge C^{-1}t^{-\frac{N}{2}\left(\frac{1}{p}-\frac{1}{q}\right)-1}
 +C^{-1}\epsilon t^{-\frac{N}{2p}}\frac{\|\nabla^2 h_0\|_{L^{q,\theta}(B(0,\delta\sqrt{t}))}}{h_0(\sqrt{t})}\\
 & \qquad\quad
 +C^{-1}\epsilon t^{-\frac{N}{2p}}\sum_{i=1}^N\frac{\|\nabla^\alpha J_{1,i}\|_{L^{q,\theta}(B(0,\delta\sqrt{t}))}}{h_1(\sqrt{t})}
 -C(N+1)\epsilon t^{-\frac{N}{2}\left(\frac{1}{p}-\frac{1}{q}\right)-1}\\
 & \ge C^{-1}t^{-\frac{N}{2}\left(\frac{1}{p}-\frac{1}{q}\right)-1}
 +C^{-1}\epsilon t^{-\frac{N}{2p}}
 \biggr[\frac{\|\nabla h_0\|_{L^{q,\theta}(B(0,\sqrt{t}))}}{h_0(\sqrt{t})}+\sum_{i=1}^N\frac{\|\nabla^\alpha J_{1,i}\|_{L^{q,\theta}(B(0,\sqrt{t}))}}{h_1(\sqrt{t})}\biggr]\\
  & \ge C^{-1}\epsilon t^{-\frac{N}{2p}}
 \biggr[\frac{\|\nabla h_0\|_{L^{q,\theta}(B(0,\sqrt{t}))}}{h_0(\sqrt{t})}+\sum_{i=1}^N\frac{\|\nabla^\alpha J_{1,i}\|_{L^{q,\theta}(B(0,\sqrt{t}))}}{h_1(\sqrt{t})}
 +t^{\frac{N}{2q}-1}\biggr]
  \ge C^{-1}\epsilon\Phi_2(t)
\end{split}
\end{equation*}
for $t\in(0,R_2^2)$. 
These imply that \eqref{eq:5.5} holds for $t\in(0,R_2^2)$. 

Combining the arguments in cases (${\mathcal A}_2$) and (${\mathcal B}_2$),  
we find $T_1>0$ such that \eqref{eq:5.5} holds for $t\in(0,T_1)$. 
Furthermore, 
we see that \eqref{eq:5.2} and \eqref{eq:5.3} hold. 
Then, for any $T_2>T_1$,  
we have 
$$
\|\nabla^\alpha e^{-tH}\|_{(L^{p,\sigma}\to L^{q,\theta})}\ge Ct^{-\frac{N}{2}\left(\frac{1}{p}-\frac{1}{q}\right)-\frac{\alpha}{2}}
\ge C^{-1}\ge C^{-1}\Phi_\alpha(t),
\quad t\in[T_1,T_2],
$$
which implies that \eqref{eq:5.5} holds for $t\in[T_1,T_2]$.

It remains to prove \eqref{eq:5.5} for $t>T_2$. 
We divide the proof in cases (${\mathcal A}_2'$) and (${\mathcal B}_2'$). 
Consider case~(${\mathcal A}_2'$). 
Similarly to \eqref{eq:5.7}, 
taking large enough $T_2$ if necessary, 
by Theorem~\ref{Theorem:4.2}~(a) and \eqref{eq:2.10} we have
\begin{equation*}
\begin{split}
 & \|\nabla^\alpha e^{-tH}\|_{(L^{p,\sigma}\to L^{q,\theta})}
\ge\|\partial_r^\alpha e^{-tH}\|_{(L^{p,\sigma}_r\to L^{q,\theta}(B(0,\sqrt{t})\setminus B(0,R_3)))}\\
 & \ge C^{-1}t^{-\frac{N}{2}}\Gamma_{p',\sigma'}(t)\frac{\|\partial_r^\alpha h_0\|_{L^{q,\theta}(B(0,\sqrt{t})\setminus B(0,R_3))}}{h_0(\sqrt{t})}\\
 & \ge C^{-1}t^{-\frac{N}{2}}\Gamma_{p',\sigma'}(t)\frac{\|\nabla^\alpha h_0\|_{L^{q,\theta}(B(0,\sqrt{t})\setminus B(0,R_3))}}{h_0(\sqrt{t})}
 \ge C^{-1}t^{-\frac{N}{2}}\Gamma_{p',\sigma'}(t)\frac{\|\nabla^\alpha h_0\|_{L^{q,\theta}(B(0,\sqrt{t}))}}{h_0(\sqrt{t})}
\end{split}
\end{equation*}
for $t\ge T_4$. 
Furthermore, by Proposition~\ref{Proposition:2.2}~(a), \eqref{eq:2.10}, and \eqref{eq:2.11} we have 
\begin{equation*}
\begin{split}
 & \sum_{i=1}^N\frac{\|\nabla^\alpha J_{1,k}\|_{L^{q,\theta}(B(0,\sqrt{t}))}}{h_1(\sqrt{t})}
 \le\sum_{i=1}^N\frac{\|\nabla^\alpha J_{1,k}\|_{L^{q,\theta}(B(0,\sqrt{t})\setminus B(0,R_3))}+C}{h_1(\sqrt{t})}\\
 & \qquad
 \le C\frac{\|\nabla^\alpha h_0\|_{L^{q,\theta}(B(0,\sqrt{t})\setminus B(0,R_3))}+1}{h_0(\sqrt{t})}
 \le C\frac{\|\nabla^\alpha h_0\|_{L^{q,\theta}(B(0,\sqrt{t}))}}{h_0(\sqrt{t})},\\
 & \frac{\|\nabla^\alpha h_0\|_{L^{q,\theta}(B(0,\sqrt{t})\setminus B(0,R_3))}}{h_0(\sqrt{t})}
 \ge C^{-1}\frac{\|h_0^{\langle\alpha\rangle}\|_{L^{q,\theta}(B(0,\sqrt{t})\setminus B(0,R_3))}}{h_0(\sqrt{t})}\\
 & \qquad
 \ge Ct^{-\frac{\alpha}{2}}\frac{\|h_0\|_{L^{q,\theta}(B(0,\sqrt{t})\setminus B(0,R_3))}}{h_0(\sqrt{t})}
 \ge Ct^{-\frac{\alpha}{2}}\Gamma_{q,\theta}(t)
\ge Ct^{\frac{N}{2q}-\frac{\alpha}{2}},
\end{split}
\end{equation*}
for $t\ge T_4$. These imply that
\eqref{eq:5.5} holds for $t\ge T_4$. 

Consider case~(${\mathcal B}_2'$). 
Similarly to \eqref{eq:5.9}, 
by Theorem~\ref{Theorem:4.1} and \eqref{eq:2.8} 
we have
\begin{equation*}
\begin{split}
 & \|\nabla^\alpha e^{-tH}\|_{(L^{p,\sigma}\to L^{q,\theta})}
 \ge \|\nabla^\alpha e^{-tH}\|_{(L^{p,\sigma}_{k,i}\to L^{q,\theta}(B(0,\delta\sqrt{t})))}\\
 & \ge t^{-\frac{N}{2}}\frac{\Gamma^k_{p',\sigma'}(t)}{h_k(\sqrt{t})}
\biggr[C^{-1}\|\nabla^\alpha J_{k,i}\|_{L^{q,\theta}(B(0,\delta\sqrt{t}))}-Ct^{-1}\|h_k^{\langle \alpha-2\rangle}\|_{L^{q,\theta}(B(0,\delta\sqrt{t}))}\biggr]_+\\
 & \ge C^{-1}t^{-\frac{N}{2p}}\frac{\|\nabla^\alpha J_{k,i}\|_{L^{q,\theta}(B(0,\sqrt{t}))}}{h_k(\sqrt{t})}
 -Ct^{-\frac{N}{2}\left(\frac{1}{p}-\frac{1}{q}\right)-\frac{\alpha}{2}}
\end{split}
\end{equation*}
for $t\ge T_4$,
where $k\in\{0,1\}$ and $i\in\{1,\dots,d_k\}$. 
Applying the same argument as in case (${\mathcal B}_2$)
we see that \eqref{eq:5.5} holds for $t\ge T_4$. 
Therefore we deduce that \eqref{eq:5.5} holds for $t>0$.  
Thus Theorem~\ref{Theorem:1.1} follows.
$\Box$
\section{Proofs of Theorems~\ref{Theorem:1.2} and \ref{Theorem:1.3}}
{\bf Proof of Theorem~\ref{Theorem:1.2}.}
Assertion~(a) follows from Theorem~\ref{Theorem:4.2}~(b). 
It suffices to prove assertion~(b). 
Assume \eqref{eq:1.9} and let $\alpha\in\{0,1,\dots,m+1\}$. 
Consider the case of $A_{1,0}<\alpha$. 
If $A_{1,0}\not\in\{0,1,2,\dots,\alpha-1\}$, 
then, by Proposition~\ref{Proposition:2.1} 
we find $R_1>0$ such that 
$$
|\nabla^\alpha h_0(|x|)|\asymp|\partial_r^\alpha h_0(|x|)|\asymp |x|^{A_{1,0}-\alpha},
\quad x\in B(0,R_1)\setminus\{0\}.
$$
This implies that 
$\partial_r^\alpha h_0\not\in L^\infty(B(0,R))$ for $R>0$. 
Then, by Theorem~\ref{Theorem:4.2}~(a) we see that 
$$
\|\nabla^\alpha e^{-tH}\|_{(L^{p,\sigma}\to L^\infty)}
\ge\|\nabla^\alpha e^{-tH}\|_{(L^{p,\sigma}_{0,1}\to L^\infty)}\ge C^{-1}\|\partial_r^\alpha e^{-tH}\|_{(L^{p,\sigma}_{0,1}\to L^\infty)}=\infty,
$$
which contradicts \eqref{eq:1.9}. This implies that 
$A_{1,0}\ge\alpha$ if $A_{1,0}\not\in\{0,1,2,\dots,\alpha-1\}$.

If $A_{1,0}\in\{1,2,\dots,\alpha-1\}$, then $\alpha\ge 1$, $\lambda_1\ge\omega_1>0$, and $0<A_{1,1}-A_{1,0}<1$ (see e.g. \cite[Lemma~4.2]{IK02}).  
Applying the above argument again, we see that $\partial_r^\alpha h_1\not\in L^\infty(B(0,R))$ for $R>0$. 
Then, by Theorem~\ref{Theorem:4.2}~(a) we see that 
$$
\|\nabla^\alpha e^{-tH}\|_{(L^{p,\sigma}\to L^\infty)}
\ge\|\nabla^\alpha e^{-tH}\|_{(L^{p,\sigma}_{1,1}\to L^\infty)}\ge C^{-1}\|\partial_r^\alpha e^{-tH_1}\|_{(L^{p,\sigma}_{0,1}\to L^\infty)}=\infty,
$$
which contradicts \eqref{eq:1.9}. This implies that $A_{1,0}\not\in\{1,2,\dots,\alpha-1\}$.
We deduce that $A_{1,0}\in[\alpha,\infty)\cup\{0\}$, that is, 
$\lambda_1\in[\omega_\alpha,\infty)\cup\{0\}$. 

Consider the case of $A_{2,0}<\alpha$. 
If $A_{2,0}\not\in\{0,1,2,\dots,\alpha-1\}$, 
then, by Proposition~\ref{Proposition:2.1} 
we find $R_2>0$ such that 
$$
|\nabla^\alpha h_0(|x|)|\asymp|\partial_r^\alpha h_0(|x|)|\asymp h_0^{\langle\alpha\rangle}(|x|),
\quad x\in B(0,R_2)^c.
$$ 
By Theorem~\ref{Theorem:4.2}~(a), Proposition~\ref{Proposition:2.2}, and \eqref{eq:1.6} we have
\begin{equation}
\label{eq:6.1}
\begin{split}
\|\nabla^\alpha e^{-tH}\|_{(L^{p,\sigma}\to L^\infty)}
 & \ge C^{-1}\|\partial_r^\alpha e^{-tH}\|_{(L^{p,\sigma}_r\to L^\infty(B(0,\delta\sqrt{t})\cap B(0,R_2)^c))}\\
 & \ge Ct^{-\frac{N}{2p}}\frac{\|h^{\langle\alpha\rangle}_0\|_{L^\infty(B(0,\delta\sqrt{t})\cap B(0,R_2)^c)}}{h_0(\sqrt{t})}
 \ge Ct^{-\frac{N}{2p}-\frac{A_{2,0}}{2}}(\log t)^{-B_0}
\end{split}
\end{equation}
for large enough $t>0$. 
This together with \eqref{eq:1.9} contradicts $A_{2,0}<\alpha$. 
Thus $A_{2,\alpha}\ge\alpha$ if $A_{2,0}\not\in\{0,1,2,\dots,\alpha-1\}$.

If $A_{2,0}\in\{1,2,\dots,\alpha-1\}$, then $\lambda_2\ge\omega_2>0$ and $0<A_{2,1}-A_{2,0}<1$ 
(see e.g. \cite[Lemma~4.2]{IK02}).
Furthermore, we find $R_3>0$ such that 
$$
|\nabla^\alpha h_1(|x|)|\asymp|\partial_r^\alpha h_1(|x|)|\asymp |x|^{A_{2,1}-\alpha},
\quad x\in B(0,R_3)^c.
$$ 
Similarly to \eqref{eq:6.1}, we have
\begin{equation*}
\begin{split}
\|\nabla^\alpha e^{-tH}\|_{(L^{p,\sigma}\to L^\infty)}
 & \ge C^{-1}\|\partial_r^\alpha e^{-tH_1}\|_{(L^{p,\sigma}_r\to L^\infty(B(0,\delta\sqrt{t})\cap B(0,R_3)^c))}\\
 & \ge Ct^{-\frac{N}{2p}}\frac{\|\partial_r^\alpha h_1\|_{L^\infty(B(0,\delta\sqrt{t})\cap B(0,R_3)^c)}}{h_1(\sqrt{t})}
 \ge Ct^{-\frac{N}{2p}-\frac{A_{2,1}}{2}}
\end{split}
\end{equation*}
for large enough $t>0$. This together with \eqref{eq:1.9} implies that $\alpha\le A_{2,1}<A_{2,0}+1\le\alpha$. 
This is a contradiction. 
So we see that $A_{2,0}\in[\alpha,\infty)\cup\{0\}$.

We prove that $A_{2,0}\in[\alpha,\infty)$ if $\alpha\ge 1$. 
Let $A_{2,0}=0$ and $R>0$. 
By Proposition~\ref{Proposition:2.2}~(a) and Theorem~\ref{Theorem:4.1} we see that 
\begin{equation}
\label{eq:6.2}
\begin{split}
 & \|\nabla e^{-tH}\|_{(L^{p,\sigma}\to L^\infty)}\ge\|\nabla e^{-tH}\|_{(L^{p,\sigma}_{0,1}\to L^\infty(B(0,R+1)\setminus B(0,R)))}\\
 & \ge t^{-\frac{N}{2}}\frac{\Gamma_{p',\sigma'}(t)}{h_0(\sqrt{t})}
 \biggr[C^{-1}\|\nabla h_0\|_{L^\infty(B(0,R+1)\setminus B(0,R))}-Ct^{-1}\|h_0^{\langle -1\rangle}\|_{L^\infty(B(0,R+1)\setminus B(0,R))}\biggr]_+\\
 & \ge C^{-1}t^{-\frac{N}{2p}}\biggr[C^{-1}\|\nabla h_0\|_{L^\infty(B(0,R+1)\setminus B(0,R))}-Ct^{-1}\biggr]_+
\end{split}
\end{equation}
for large enough $t>0$. This together with \eqref{eq:1.9} implies that $\|\nabla h_0\|_{L^\infty(B(0,R+1)\setminus B(0,R))}=0$. 
Since $R$ is arbitrary, we observe that $h_0$ is a constant function in ${\bf R}^N$, 
which contradicts that $V\not\equiv 0$ in ${\bf R}^N$ (see \eqref{eq:1.5}). Thus $A_{2,0}\not=0$. 
Therefore we see that $A_{2,0}\in[\alpha,\infty)$, that is, $H$ is subcritical and $\lambda_2\ge\omega_\alpha$ 
(see \eqref{eq:1.3}). 
Thus Theorem~\ref{Theorem:1.2} follows.
$\Box$
\vspace{5pt}
\newline
{\bf Proof of Theorem~\ref{Theorem:1.3}.}
Since $\alpha$ is arbitrary, 
we apply a similar argument as in that of the proof of Theorem~\ref{Theorem:1.2} to see that $\lambda_1=\lambda_2=0$. 
Furthermore, for $R>0$, we have
\begin{equation*}
\|\nabla e^{-tH}\|_{(L^{p,\sigma}\to L^{q,\theta})}
\ge C^{-1}t^{-\frac{N}{2p}}\biggr[C^{-1}\|\nabla h_0\|_{L^{q,\theta}(B(0,R+1)\setminus B(0,R))}-Ct^{-1}\biggr]_+
\end{equation*}
for large enough $t>0$, instead of \eqref{eq:6.2}. 
Then, similarly to the proof of Theorem~\ref{Theorem:1.2}, 
we observe from \eqref{eq:1.10} that $h_0$ is a constant function in ${\bf R}^N$. 
This means that $V\equiv 0$ in ${\bf R}^N$. Thus Theorem~\ref{Theorem:1.3} follows.
$\Box$
\section{Applications}
We study the decay of $\|\nabla^\alpha e^{-tH}\|_{(L^{p,\sigma}\to L^{q,\theta})}$ 
for some typical potentials. 
\subsection{Hardy potentials}
In this subsection we consider the case when
\begin{equation}
\label{eq:7.1}
V(r)=\lambda r^{-2},\quad r>0,
\end{equation}
with $\lambda\ge\lambda_*=-(N-2)^2/4$ and $\lambda\not=0$. 
(See Remark~\ref{Remark:3.1} for the case of $\lambda=0$.) 
Set 
$$
A_k:=\frac{-(N-2)+\sqrt{(N-2)^2+4(\lambda+\omega_k)}}{2},\quad k\in\{0,1,2,\dots\}.
$$
We remark that, under \eqref{eq:7.1}, 
$H=-\Delta+V$ is subcritical if $\lambda>\lambda_*$ and it is critical if $\lambda=\lambda_*$. 
Then 
$$
A_{1,k}=A_{2,k}=A_k,
\qquad
h_k(|x|)=|x|^{A_k},
$$
for $k\in\{0,1,2,\dots\}$ and $\lambda\ge\lambda_*$. Furthermore, $A_0\not=0$ by $\lambda\not=0$. 
\begin{theorem}
\label{Theorem:7.1} 
Assume \eqref{eq:7.1} with $\lambda\ge\lambda_*=-(N-2)^2/4$ and $\lambda\not=0$. 
Let $\alpha\in\{0,1,\dots\}$ and $(p,q,\sigma,\theta)\in\Lambda$. 
For $A\in{\bf R}$, set  $f_A(x):=|x|^A$ for $\in{\bf R}^N\setminus\{0\}$.
\begin{itemize}
  \item[{\rm (a)}] 
  There exists $C_1>0$ such that 
  $$
  \|\nabla^\alpha e^{-tH}\|_{(L^{p,\sigma}\to L^{q,\theta})}
  \ge C_1^{-1}t^{-\frac{N}{2}\left(\frac{1}{p}-\frac{1}{q}\right)-\frac{\alpha}{2}},\quad t>0. 
  $$
  \item[{\rm (b)}]
  Assume that $A_0\not\in\{2,4,\dots\}$.  
  Then $\|\nabla^\alpha e^{-tH}\|_{(L^{p,\sigma}\to L^{q,\theta})}<\infty$ for some $t>0$ if and only if 
  \begin{equation}
  \label{eq:7.2}
  f_{A_0}\in L^{p',\sigma'}(B(0,1)),\qquad f_{A_0}^{\langle\alpha\rangle}\in L^{q,\theta}(B(0,1)).
  \end{equation}
  Furthermore, under assumption~\eqref{eq:7.2}, there exists $C_2>0$ that
  \begin{equation}
  \label{eq:7.3}
  \|\nabla^\alpha e^{-tH}\|_{(L^{p,\sigma}\to L^{q,\theta})}
  \le C_2t^{-\frac{N}{2}\left(\frac{1}{p}-\frac{1}{q}\right)-\frac{\alpha}{2}},\quad t>0.
  \end{equation} 
  \item[{\rm (c)}]
  Assume that $A_0=2\gamma$ for some $\gamma\in\{1,2,\dots\}$. 
  Then $\|\nabla^\alpha e^{-tH}\|_{(L^{p,\sigma}\to L^{q,\theta})}<\infty$ for some $t>0$ if and only if 
  \begin{equation}
  \label{eq:7.4}
  \mbox{either}\quad \alpha\le 2\gamma\quad\mbox{or}\quad f_{A_1}^{\langle\alpha\rangle}\in L^{q,\theta}(B(0,1)).
  \end{equation}
  Furthermore, under assumption~\eqref{eq:7.4},  
  inequality~\eqref{eq:7.3} holds  for some $C_2>0$.
\end{itemize}
\end{theorem}
{\bf Proof.} 
Assertion~(a) follows from Theorem~\ref{Theorem:4.2}~(b).
We prove assertions~(b) and (c). 
Let $(k,i)\in{\mathcal K}$. 
By Theorem~\ref{Theorem:4.1} we have
\begin{equation}
\label{eq:7.5}
\begin{split}
 \|\nabla^\alpha e^{-tH}\|_{(L^{p,\sigma}\to L^{q,\theta})}
 & \ge \|\nabla^\alpha e^{-tH}\|_{(L^{p,\sigma}_{k,i}\to L^{q,\theta})}\\
 & \ge t^{-\frac{N}{2}}\frac{\Gamma_{p',\sigma'}(t)}{h_k(\sqrt{t})}
 \biggr[C^{-1}\|\nabla^\alpha J_{k,i}\|_{L^{q,\theta}(E)}-Ct^{-1}\|h_k^{\langle \alpha-2\rangle}\|_{L^{q,\theta}(E)}\biggr]_+
\end{split}
\end{equation}
for measurable sets $E\subset B(0,\delta\sqrt{t})$ and $t>0$, where $\delta\in(0,1)$ is as in Theorem~\ref{Theorem:4.1}. 
Furthermore, 
for any $(k,i)\in{\mathcal K}$ and $\alpha, n\in\{0,1,2,\dots\}$, by Lemma~\ref{Lemma:2.1}~(b) we have
\begin{equation}
\label{eq:7.6}
t^{-n}\frac{|\nabla^\alpha J_{k,i}^n(x)|}{h_k(\sqrt{t})}\le Ct^{-n}|x|^{2n-\alpha} \frac{h_k(|x|)}{h_k(\sqrt{t})}
\le Ct^{-\frac{A_k}{2}-n}|x|^{A_k+2n-\alpha},
\quad x\in{\bf R}^N\setminus\{0\},\,\,t>0.
\end{equation}

Assume that $A_0\not\in\{2,4,\dots\}$. 
Since $|\nabla^\alpha |x|^{A_0}|\not\equiv 0$ in ${\bf R}^N\setminus\{0\}$, 
taking small enough $\epsilon\in(0,\delta)$, 
by \eqref{eq:7.5} we have 
\begin{equation}
\label{eq:7.7}
\begin{split}
\|\nabla^\alpha e^{-tH}\|_{(L^{p,\sigma}\to L^{q,\theta})}
 & \ge t^{-\frac{N}{2}}\frac{\Gamma_{p',\sigma'}(t)}{h_0(\sqrt{t})}
\biggr[C^{-1}\|\nabla^\alpha f_{A_0}\|_{L^{q,\theta}(D_R(t))}-C(\epsilon\delta)^2\|f_{A_0}^{\langle \alpha\rangle}\|_{L^{q,\theta}(D_R(t))}\biggr]_+\\
 & \ge C^{-1}t^{-\frac{N}{2}}\frac{\Gamma_{p',\sigma'}(t)}{h_0(\sqrt{t})}\|f^{\langle\alpha\rangle}_{A_0}\|_{L^{q,\theta}(D_R(t))}
\end{split}
\end{equation}
for small enough $R>0$, where $D_R(t):=B(0,\epsilon\delta\sqrt{t})\cap B(0,R)^c$. 
Since $R$ is arbitrary, 
we see that \eqref{eq:7.2} holds if 
$\|\nabla^\alpha e^{-tH}\|_{(L^{p,\sigma}\to L^{q,\theta})}<\infty$ for some $t>0$.
Furthermore, 
since $A_k+2n-\alpha\ge A_0-\alpha$, 
if $f_{A_0}^{\langle\alpha\rangle}\in L^{q,\theta}(B(0,1))$, by \eqref{eq:2.4} we see that 
$$
\sum_{0\le k+2n\le\alpha}\sum_{i=1}^{d_k}
 t^{-n}\frac{\|\nabla^\alpha J_{k,i}^n\|_{L^{q,\theta}(B(0,\sqrt{t}))}}{h_k(\sqrt{t})}
\le Ct^{-\frac{A_0}{2}}\|f_{A_0}^{\langle\alpha\rangle}\|_{L^{q,\theta}(B(0,\sqrt{t}))}
\le Ct^{\frac{N}{2q}-\frac{\alpha}{2}},\quad t>0.
$$
This together with Theorem~\ref{Theorem:3.1}~(b) implies assertion~(b). 
Similarly, we see that assertion~(c) holds in the case of $\alpha\le2\gamma$.

Consider the case when $A_0=2\gamma$ with $\gamma\in\{1,2,\dots\}$ and $\alpha>2\gamma\ge 2$. 
Assume that $f_{A_1}^{\langle\alpha\rangle}\in L^{q,\theta}(B(0,1))$.
It follows from $\lambda>0$ that $A_1<A_0+1$ (see e.g. \cite[Lemma~4.2]{IK02}).  
Since $\nabla^\alpha |x|^{A_1}\not\equiv 0$ in ${\bf R}^N\setminus\{0\}$, 
similarly to \eqref{eq:7.7}, taking small enough $\epsilon\in(0,\delta)$ if necessary, we have 
\begin{equation*}
\begin{split}
\|\nabla^\alpha e^{-tH}\|_{(L^{p,\sigma}\to L^{q,\theta})}
 & \ge t^{-\frac{N}{2}}\frac{\Gamma_{p',\sigma'}(t)}{h_1(\sqrt{t})}
\biggr[C^{-1}\|\nabla^\alpha f_{A_1}\|_{L^{q,\theta}(D_R(t))}-C(\epsilon\delta)^2\|f_{A_1}^{\langle \alpha\rangle}\|_{L^{q,\theta}(D_R(t))}\biggr]_+\\
 & \ge C^{-1}t^{-\frac{N}{2}}\frac{\Gamma_{p',\sigma'}(t)}{h_1(\sqrt{t})}\|f^{\langle\alpha\rangle}_{A_1}\|_{L^{q,\theta}(D_R(t))}
\end{split}
\end{equation*}
for small enough $R>0$. 
This implies that 
$f_{A_1}^{\langle\alpha\rangle}\in L^{q,\theta}(B(0,1))$ if 
$\|\nabla^\alpha e^{-tH}\|_{(L^{p,\sigma}\to L^{q,\theta})}<\infty$ for some $t>0$. 
Furthermore, since $A_1<A_0+1$, 
if $f_{A_1}^{\langle\alpha\rangle}\in L^{q,\theta}(B(0,1))$, 
then $f_{A_0}^{\langle\alpha-2\rangle}\in L^{q,\theta}(B(0,1))$. 
On the other hand, since $h_{A_0}$ is a homogeneous polynomial of degree $2\gamma$, 
we see that $|\nabla^\alpha J_{0,i}^0(x)|=C|\nabla^\alpha f_{A_0}(|x|)|\equiv 0$ in ${\bf R}^N$ if $\alpha>2\gamma$. 
These together with \eqref{eq:7.6} imply that 
\begin{equation*}
\begin{split}
 & \sum_{0\le k+2n\le\alpha}\sum_{i=1}^{d_k}
 t^{-n}\frac{\|\nabla^\alpha J_{k,i}^n\|_{L^{q,\theta}(B(0,\sqrt{t}))}}{h_k(\sqrt{t})}\\
 & =\sum_{0<2n\le\alpha}
 t^{-n}\frac{\|\nabla^\alpha J_{0,i}^n\|_{L^{q,\theta}(B(0,\sqrt{t}))}}{h_0(\sqrt{t})}
 +\sum_{\substack{0\le k+2n\le\alpha\\k\ge 1}}\sum_{i=1}^{d_k}
 t^{-n}\frac{\|\nabla^\alpha J_{k,i}^n\|_{L^{q,\theta}(B(0,\sqrt{t}))}}{h_k(\sqrt{t})}\\
 & \le Ct^{-\frac{A_0}{2}-1}\|f_{A_0}^{\langle\alpha-2\rangle}\|_{L^{q,\theta}(B(0,\sqrt{t}))}
 +Ct^{-\frac{A_1}{2}}\|f_{A_1}^{\langle\alpha\rangle}\|_{L^{q,\theta}(B(0,\sqrt{t}))}
 \le Ct^{\frac{N}{2q}-\frac{\alpha}{2}},\quad t>0.
\end{split}
\end{equation*}
Then, by Theorem~\ref{Theorem:3.1}~(b) we have assertion~(c). Thus Theorem~\ref{Theorem:7.1} follows. 
$\Box$
\subsection{Bounded potentials}
We consider the case when $V\in C^m([0,\infty))$. Then $\lambda_1=0$ and $h_k$ has no singularity at $x=0$. 
In Theorem~\ref{Theorem:7.2} we treat the following two cases:
$$
{\rm (A)}\quad\mbox{either $A_{2,0}\not\in\{0,2,4,\dots\}$ or $A_{2,0}\ge\alpha$},
\qquad
{\rm (B)}\quad\mbox{$A_{2,0}\in\{2,4,\dots\}$ and $A_{2,0}<\alpha$}.
$$
The case when $A_{2,0}=0$ and $\alpha\ge 1$ is discussed later.
\begin{theorem}
\label{Theorem:7.2}
Let $V\in C^m([0,\infty))$ for some $m\in\{0,1,2,\dots\}$. 
Assume conditions~{\rm ($\mbox{V}_m$)} and~{\rm (N')}. 
Let $\alpha\in\{0,1,\dots,m+1\}$, $(p,q,\sigma,\theta)\in\Lambda$, and $T>0$. 
Then 
\begin{equation}
\label{eq:7.8}
\|\nabla^\alpha e^{-tH}\|_{(L^{p,\sigma}\to L^{q,\theta})}
\asymp t^{-\frac{N}{2}\left(\frac{1}{p}-\frac{1}{q}\right)-\frac{\alpha}{2}},
\quad 0<t\le T.
\end{equation}
Furthermore, 
\begin{equation}
\label{eq:7.9}
\begin{split}
 & \|\nabla^\alpha e^{-tH}\|_{(L^{p,\sigma}\to L^{q,\theta})}\\
 & \asymp
\left\{
\begin{array}{ll}
\displaystyle{t^{-\frac{N}{2}}\Gamma_{p',\sigma'}(t)\frac{\|\tilde{h}_0^{\langle\alpha\rangle}\|_{L^{q,\theta}(B(0,\sqrt{t}))}}{h_0(\sqrt{t})}}
 & \mbox{in case~{\rm (A)}},\vspace{5pt}\\
\displaystyle{t^{-\frac{N}{2}}\Gamma_{p',\sigma'}(t)
\biggr[\frac{\|\nabla^\alpha h_0\|_{L^{q,\theta}(B(0,\sqrt{t}))}}{h_0(\sqrt{t})}
+\frac{\|\tilde{h}_1^{\langle\alpha\rangle}\|_{L^{q,\theta}(B(0,\sqrt{t}))}}{h_1(\sqrt{t})}}\biggr]
 & \mbox{in case {\rm (B)}},
\end{array}
\right.
\end{split}
\end{equation}
for $t>T$. Here $\tilde{h}_k^{\langle\alpha\rangle}(x)=(1+|x|)^{-\alpha}h_k(|x|)$ for $k=0,1$.
\end{theorem}
{\bf Proof.}
It follows from $V\in C^m([0,\infty))$ that, 
for any $R>0$ and $k\in\{0,1,2,\dots\}$, 
$|\nabla^\alpha h_k|$ is bounded in $B(0,R)$ for $\alpha\in\{0,1,\dots,m+1\}$. 
Then, by Theorems~\ref{Theorem:3.1} and \ref{Theorem:4.2}~(b) 
we easily obtain relation~\eqref{eq:7.8}.
In case~(A) we find $R_1>0$ such that 
$$
|\nabla^\alpha h_0(|x|)|\asymp h^{\langle\alpha\rangle}_0(|x|),\quad x\in B(0,R_1)^c.
$$
Furthermore, Proposition~\ref{Proposition:2.2}~(a) implies that 
$$
\frac{\|\tilde{h}_0^{\langle\alpha\rangle}\|_{L^{q,\theta}(B(0,\sqrt{t}))}}{h_0(\sqrt{t})}
\ge (1+\sqrt{t})^{-\alpha}\frac{\|h_0\|_{L^{q,\theta}(B(0,\sqrt{t}))}}{h_0(\sqrt{t})}
\ge Ct^{\frac{N}{2q}-\frac{\alpha}{2}},\quad t\ge T.
$$
Since $\nabla^\alpha h_k$ is bounded in $B(0,R_1)$, 
by Lemma~\ref{Lemma:2.1}~(b) we apply Theorem~\ref{Theorem:3.1}, Theorem~\ref{Theorem:4.2}, and Proposition~\ref{Proposition:2.2}~(c) 
to obtain
\begin{equation}
\label{eq:7.10}
\begin{split}
\|\nabla^\alpha e^{-tH}\|_{(L^{p,\sigma}\to L^{q,\theta})}
 & \asymp t^{-\frac{N}{2}}\Gamma_{p',\sigma'}(t)
\biggr[\frac{\|\tilde{h}_0^{\langle\alpha\rangle}\|_{L^{q,\theta}(B(0,\sqrt{t}))}}{h_0(\sqrt{t})}+t^{\frac{N}{2q}-\frac{\alpha}{2}}\biggr]\\
 & \asymp t^{-\frac{N}{2}}\Gamma_{p',\sigma'}(t)
\frac{\|\tilde{h}_0^{\langle\alpha\rangle}\|_{L^{q,\theta}(B(0,\sqrt{t}))}}{h_0(\sqrt{t})},
\quad t>T.
\end{split}
\end{equation}
On the other hand, 
in case~(B) we see that $A_{2,0}\in\{2,4,\dots\}$ and $0<A_{2,1}-A_{2,0}<1$. 
Then we find $R_2>0$ such that
$$
|\nabla^\alpha h_1(x|)|\asymp h^{\langle\alpha\rangle}_1(|x|),\quad x\in B(0,R_2)^c.
$$
Then, by Lemma~\ref{Lemma:2.1} and Proposition~\ref{Proposition:2.2}~(c) we have
\begin{equation*}
\begin{split}
t^{-n}\frac{|\nabla^\alpha J_{k,i}^n(x)|}{h_k(\sqrt{t})} & \le C\frac{h_1^{\langle\alpha\rangle}(|x|)}{h_1(\sqrt{t})}
\quad\mbox{if}\quad k\ge 1,\\
t^{-n}\frac{|\nabla^\alpha J_{0,i}^n(x)|}{h_0(\sqrt{t})} & \le C\frac{|x|^2h_0^{\langle\alpha\rangle}(|x|)}{th_0(\sqrt{t})}
\le C\frac{h_1^{\langle\alpha\rangle}(|x|)}{h_1(\sqrt{t})}
 \quad\mbox{if}\quad n\ge 1,
\end{split}
\end{equation*}
for $x\in B(0,\sqrt{t})\cap B(0,R_2)^c$ and $0\le k+2n\le\alpha$. 
Since 
$$
\frac{\|\tilde{h}_1^{\langle\alpha\rangle}\|_{L^{q,\theta}(B(0,\sqrt{t}))}}{h_1(\sqrt{t})}\ge Ct^{\frac{N}{2q}-\frac{\alpha}{2}},\quad t\ge T,
$$
similarly to \eqref{eq:7.10}, we apply Theorem~\ref{Theorem:3.1}, Theorem~\ref{Theorem:4.2}, and Proposition~\ref{Proposition:2.2}~(c) 
to obtain inequality~\eqref{eq:7.9} in case~(B). 
Thus Theorem~\ref{Theorem:7.2} follows.
$\Box$
\vspace{5pt}
\newline
Under the assumptions of Theorem~\ref{Theorem:7.2}, 
the exact large time decay rate of $\|\nabla^\alpha e^{-tH}\|_{(L^{p,\sigma}\to L^{q,\theta})}$ 
for bounded potentials is obtained by the combination of \eqref{eq:7.8}, \eqref{eq:7.9}, and 
\begin{equation*}
\Gamma_{p',\sigma'}(t)\asymp
\left\{
\begin{array}{ll}
t^{-\frac{A_{2,0}}{2}}(\log t)^{-B_0}\quad & \mbox{if $p<p_*$},\\
t^{-\frac{A_{2,0}}{2}}(\log t)^{\frac{1}{\sigma'}}\quad & \mbox{if $p=p_*$},\\
t^{\frac{N}{2}\left(1-\frac{1}{p}\right)}\quad & \mbox{if $p>p_*$},
\end{array}
\right.
\frac{\|\tilde{h}_0^{\langle\alpha\rangle}\|_{L^{q,\theta}(B(0,\sqrt{t}))}}{h_0(\sqrt{t})}
\asymp
\left\{
\begin{array}{ll}
t^{\frac{N}{2q}-\frac{\alpha}{2}}\quad & \mbox{if $q<q_\alpha$},\\
t^{-\frac{A}{2}}(\log t)^{\frac{1}{\theta}}\quad & \mbox{if $q=q_\alpha$},\\
t^{-\frac{A}{2}}(\log t)^{-B_0}\quad & \mbox{if $q>q_\alpha$},
\end{array}
\right.
\end{equation*}
for large enough $t>0$. 
Here $p_*:=N/(N+A_{2,0})$, $q_\alpha:=N/(\alpha-A_{2,0})$ and $B_0$ is as in \eqref{eq:1.3}. 

Finally, we discuss the large time decay of $\|\nabla^\alpha e^{-tH}\|_{(L^{p,\sigma}\to L^{q,\theta})}$  
in the case when $\lambda_2=0$ and $\alpha\in\{1,2,\dots,m+1\}$ under some additional conditions.
\begin{theorem}
\label{Theorem:7.3}
Let $m\in\{0,1,2,\dots\}$ and $V\in C^m([0,\infty))$ satisfy 
\begin{equation}
\label{eq:7.11}
\frac{d^\ell}{dr^\ell}V(r)=a(1+o(1))\frac{d^\ell}{dr^\ell}r^{-\kappa}\quad\mbox{as}\quad r\to\infty
\end{equation}
for $\ell\in\{0,1,\dots,m\}$, where $a\not=0$ and $\kappa>2$.
Assume conditions~{\rm ($\mbox{V}_m$)} and~{\rm (N')}. 
Let $(p,q,\sigma,\theta)\in\Lambda$ and $\alpha\in\{1,\dots,m+1\}$. 
Set 
$$
\eta_\alpha(|x|):=
\left\{
\begin{array}{ll}
(1+|x|)^{-\kappa+2-\alpha} & \mbox{if $2<\kappa<N$},\\
(1+|x|)^{-N+2-\alpha}\log(|x|+2) & \mbox{if $\kappa=N$},\\
(1+|x|)^{-N+2-\alpha} & \mbox{if $\kappa>N$}.
\end{array}
\right.
$$
Then 
$$
\|\nabla^\alpha e^{-tH}\|_{(L^{p,\sigma}\to L^{q,\theta})}\asymp
t^{-\frac{N}{2p}}\biggr[\|\eta_\alpha\|_{L^{q,\theta}(B(0,\sqrt{t}))}+t^{\frac{N}{2q}-\frac{\alpha}{2}}\biggr]
$$
for large enough $t>0$. 
\end{theorem}
{\bf Proof.}
For any $k\in\{0,1,2,\dots\}$, It follows from $A_{2,k}=k$ that 
$h_k(r)=F_k(r)$ for $r>0$,
where
$$
F_k(r):=r^k\left[1+\int_0^r s^{-2k-N+1}\biggr(\int_0^s \tau^{k+N-1}V(\tau)h_k(\tau)\,d\tau\biggr)\,ds\right]. 
$$
Indeed, $F_k$ satisfies \eqref{eq:1.5}. 
Since $h_k\in C^\alpha({\bf R}^N)$ and $|x|^kQ_{k,i}(x/|x|)$ is a homogeneous polynomial of degree $k$, 
by \eqref{eq:7.11} we see that 
\begin{equation*}
\begin{split}
 & \|\nabla^\alpha h_0\|_{L^{q,\theta}(B(0,\sqrt{t}))}\asymp \|\eta_\alpha\|_{L^{q,\theta}(B(0,\sqrt{t}))},\\
 & t^{-n}\frac{\|\nabla^\alpha J_{k,i}^n\|_{L^{q,\theta}(B(0,\sqrt{t}))}}{h_k(\sqrt{t})}
\le C\|\eta_\alpha\|_{L^{q,\theta}(B(0,\sqrt{t}))}+Ct^{\frac{N}{2q}},
\end{split}
\end{equation*}
for large enough $t>0$, where $(k,i)\in{\mathcal K}$ and $n\in\{0,1,\dots\}$ with $0\le k+2n\le\alpha$. 
These together with Theorems~\ref{Theorem:3.1} and \ref{Theorem:4.2} imply 
the desired inequality. Thus Theorem~\ref{Theorem:7.3} follows.
$\Box$
\vspace{5pt}
\newline
Similarly we have:
\begin{theorem}
\label{Theorem:7.4}
Let $V\in C^m([0,\infty))$ for some $m\in\{0,1,2,\dots\}$. 
Assume conditions~{\rm ($\mbox{V}_m$)} and $A_{2,0}=0$.  
Furthermore, assume that 
$$
\tau^{N-1}V(\tau)\in L^1((0,\infty)),\qquad
\int_0^\infty \tau^{N-1}V(\tau)h_0(\tau)\,d\tau\not=0.
$$
Let $(p,q,\sigma,\theta)\in\Lambda$ and $\alpha\in\{1,2,\dots\}$. Then 
$$
\|\nabla^\alpha e^{-tH}\|_{(L^{p,\sigma}\to L^{q,\theta})}\asymp
t^{-\frac{N}{2p}}\biggr[\|\tilde{\eta}_\alpha\|_{L^{q,\theta}(B(0,\sqrt{t}))}+t^{\frac{N}{2q}-\frac{\alpha}{2}}\biggr]
$$
for large enough $t>0$, where $\tilde{\eta}_\alpha(|x|):=(1+|x|)^{-N+2-\alpha}$.
\end{theorem}
By Theorems~\ref{Theorem:7.3} and \ref{Theorem:7.4} we see that 
the large time decay rate of $\|\nabla^\alpha e^{t\Delta}\|_{(L^{p,\sigma}\to L^{q,\theta})}$ varies discontinuously
with respect to perturbations of the potential~$V$ if $\alpha\not=0$. 
See also Theorem~\ref{Theorem:1.3}. 
\vspace{5pt}

\noindent
{\bf Acknowledgements.} 
The first author was supported in part by the Grant-in-Aid for Scientific Research (S)(No.~19H05599)
from Japan Society for the Promotion of Science. 
The second author was supported in part by the Grant-in-Aid for JSPS Fellows (No.~20J10379). 
  

\begin{bibdiv}
\begin{biblist}
\bib{LL}{article}{
   author={Angiuli, Luciana},
   author={Lorenzi, Luca},
   title={On the estimates of the derivatives of solutions to nonautonomous
   Kolmogorov equations and their consequences},
   journal={Riv. Math. Univ. Parma (N.S.)},
   volume={7},
   date={2016},
   pages={421--471},
}
\bib{Bar}{article}{
   author={Barbatis, G.},
   author={Filippas, S.},
   author={Tertikas, A.},
   title={Critical heat kernel estimates for Schr\"{o}dinger operators via
   Hardy-Sobolev inequalities},
   journal={J. Funct. Anal.},
   volume={208},
   date={2004},
   pages={1--30},
}
\bib{BS}{book}{
   author={Bennett, Colin},
   author={Sharpley, Robert},
   title={Interpolation of operators},
   series={Pure and Applied Mathematics},
   volume={129},
   publisher={Academic Press, Inc., Boston, MA},
   date={1988},
   pages={xiv+469},
}
\bib{BFL}{article}{
   author={Bertoldi, Marcello},
   author={Fornaro, Simona},
   author={Lorenzi, Luca},
   title={Gradient estimates for parabolic problems with unbounded
   coefficients in non convex unbounded domains},
   journal={Forum Math.},
   volume={19},
   date={2007},
   pages={603--632},
}
\bib{CK}{article}{
   author={Chavel, Isaac},
   author={Karp, Leon},
   title={Large time behavior of the heat kernel: the parabolic
   $\lambda$-potential alternative},
   journal={Comment. Math. Helv.},
   volume={66},
   date={1991},
   pages={541--556},
}
\bib{CUR}{article}{
   author={Cruz-Uribe, D.},
   author={Rios, Cristian},
   title={Gaussian bounds for degenerate parabolic equations},
   journal={J. Funct. Anal.},
   volume={255},
   date={2008},
   pages={283--312; {\it Corrigendum} in J. Funct. Anal. {\bf 267} (2014), 3507--3513},
}
\bib{Dav}{book}{
   author={Davies, E. B.},
   title={Heat kernels and spectral theory},
   series={Cambridge Tracts in Mathematics},
   volume={92},
   publisher={Cambridge University Press, Cambridge},
   date={1989},
   pages={x+197},
}
\bib{DS}{article}{
   author={Davies, E. B.},
   author={Simon, B.},
   title={$L^p$ norms of noncritical Schr\"{o}dinger semigroups},
   journal={J. Funct. Anal.},
   volume={102},
   date={1991},
   pages={95--115},
}
\bib{Gra}{book}{
   author={Grafakos, Loukas},
   title={Classical Fourier analysis},
   series={Graduate Texts in Mathematics},
   volume={249},
   edition={3},
   publisher={Springer, New York},
   date={2014},
   pages={xviii+638},
}
\bib{Gri}{book}{
   author={Grigor'yan, Alexander},
   title={Heat kernel and analysis on manifolds},
   series={AMS/IP Studies in Advanced Mathematics},
   volume={47},
   publisher={American Mathematical Society, Providence, RI; International
   Press, Boston, MA},
   date={2009},
   pages={xviii+482},
}
\bib{GS}{article}{
   author={Grigor\cprime yan, Alexander},
   author={Saloff-Coste, Laurent},
   title={Dirichlet heat kernel in the exterior of a compact set},
   journal={Comm. Pure Appl. Math.},
   volume={55},
   date={2002},
   pages={93--133},
}
\bib{GS2}{article}{
   author={Grigor'yan, Alexander},
   author={Saloff-Coste, Laurent},
   title={Stability results for Harnack inequalities},
   journal={Ann. Inst. Fourier (Grenoble)},
   volume={55},
   date={2005},
   pages={825--890},
}
\bib{IIY01}{article}{
   author={Ioku, Norisuke},
   author={Ishige, Kazuhiro},
   author={Yanagida, Eiji},
   title={Sharp decay estimates of $L^q$-norms for nonnegative Schr\"{o}dinger
   heat semigroups},
   journal={J. Funct. Anal.},
   volume={264},
   date={2013},
   pages={2764--2783},
}
\bib{IIY02}{article}{
   author={Ioku, Norisuke},
   author={Ishige, Kazuhiro},
   author={Yanagida, Eiji},
   title={Sharp decay estimates in Lorentz spaces for nonnegative
   Schr\"{o}dinger heat semigroups},
   journal={J. Math. Pures Appl. (9)},
   volume={103},
   date={2015},
   pages={900--923},
}
\bib{Ishige}{article}{
   author={Ishige, Kazuhiro},
   title={Gradient estimates for the heat equation in the exterior domains
   under the Neumann boundary condition},
   journal={Differential Integral Equations},
   volume={22},
   date={2009},
   pages={401--410},
}
\bib{IK01}{article}{
   author={Ishige, Kazuhiro},
   author={Kabeya, Yoshitsugu},
   title={Decay rates of the derivatives of the solutions of the heat
   equations in the exterior domain of a ball},
   journal={J. Math. Soc. Japan},
   volume={59},
   date={2007},
   pages={861--898},
}
\bib{IK02}{article}{
   author={Ishige, Kazuhiro},
   author={Kabeya, Yoshitsugu},
   title={Large time behaviors of hot spots for the heat equation with a
   potential},
   journal={J. Differential Equations},
   volume={244},
   date={2008},
   pages={2934--2962; {\it Corrigendum} in J. Differential Equations {\bf 245} (2008), 2352--2354},
}
\bib{IK05}{article}{
   author={Ishige, Kazuhiro},
   author={Kabeya, Yoshitsugu},
   title={$L^p$ norms of nonnegative Schr\"{o}dinger heat semigroup and the
   large time behavior of hot spots},
   journal={J. Funct. Anal.},
   volume={262},
   date={2012},
   pages={2695--2733},
}
\bib{IKM}{article}{
   author={Ishige, Kazuhiro},
   author={Kabeya, Yoshitsugu},
   author={Mukai, Asato},
   title={Hot spots of solutions to the heat equation with inverse square
   potential},
   journal={Appl. Anal.},
   volume={98},
   date={2019},
   pages={1843--1861},
}
\bib{IKO}{article}{
   author={Ishige, Kazuhiro},
   author={Kabeya, Yoshitsugu},
   author={Ouhabaz, El Maati},
   title={The heat kernel of a Schr\"{o}dinger operator with inverse square
   potential},
   journal={Proc. Lond. Math. Soc. (3)},
   volume={115},
   date={2017},
   pages={381--410},
}
\bib{IM}{article}{
   author={Ishige, Kazuhiro},
   author={Mukai, Asato},
   title={Large time behavior of solutions of the heat equation with inverse
   square potential},
   journal={Discrete Contin. Dyn. Syst.},
   volume={38},
   date={2018},
   pages={4041--4069},
}
\bib{IT01}{article}{
   author={Ishige, Kazuhiro},
   author={Tateishi, Yujiro},
   title={Decay estimates for Schr\"odinger heat semigroup with inverse square potential in Lorentz spaces},
   journal={preprint (arXiv:2009.07001)},
}
\bib{LSU}{book}{
   author={Lady\v{z}enskaja, O. A.},
   author={Solonnikov, V. A.},
   author={Ural\cprime ceva, N. N.},
   title={Linear and quasilinear equations of parabolic type},
   series={Translated from the Russian by S. Smith. Translations of
   Mathematical Monographs, Vol. 23},
   publisher={American Mathematical Society, Providence, R.I.},
   date={1968},
   pages={xi+648},
}
\bib{LY}{article}{
   author={Li, Peter},
   author={Yau, Shing-Tung},
   title={On the parabolic kernel of the Schr\"{o}dinger operator},
   journal={Acta Math.},
   volume={156},
   date={1986},
   pages={153--201},
}
\bib{Lieb}{book}{
   author={Lieberman, Gary M.},
   title={Second order parabolic differential equations},
   publisher={World Scientific Publishing Co., Inc., River Edge, NJ},
   date={1996},
   pages={xii+439},
}
\bib{LS}{article}{
   author={Liskevich, Vitali},
   author={Sobol, Zeev},
   title={Estimates of integral kernels for semigroups associated with
   second-order elliptic operators with singular coefficients},
   journal={Potential Anal.},
   volume={18},
   date={2003},
   pages={359--390},
}
\bib{MR2064932}{article}{
   author={Milman, Pierre D.},
   author={Semenov, Yu. A.},
   title={Global heat kernel bounds via desingularizing weights},
   journal={J. Funct. Anal.},
   volume={212},
   date={2004},
   pages={373--398; {\it Corrigendum} in J. Funct. Anal. {\bf 220} (2005), 238--239},
}
\bib{MT1}{article}{
   author={Moschini, Luisa},
   author={Tesei, Alberto},
   title={Harnack inequality and heat kernel estimates for the Schr\"{o}dinger
   operator with Hardy potential},
   journal={Atti Accad. Naz. Lincei Cl. Sci. Fis. Mat. Natur. Rend. Lincei
   (9) Mat. Appl.},
   volume={16},
   date={2005},
   pages={171--180 (2006)},
}
\bib{MT2}{article}{
   author={Moschini, Luisa},
   author={Tesei, Alberto},
   title={Parabolic Harnack inequality for the heat equation with
   inverse-square potential},
   journal={Forum Math.},
   volume={19},
   date={2007},
   pages={407--427},
}
\bib{M0}{article}{
   author={Murata, Minoru},
   title={Positive solutions and large time behaviors of Schr\"{o}dinger
   semigroups, Simon's problem},
   journal={J. Funct. Anal.},
   volume={56},
   date={1984},
   pages={300--310},
}
\bib{M}{article}{
   author={Murata, Minoru},
   title={Structure of positive solutions to $(-\Delta+V)u=0$ in ${\bf
   R}^n$},
   journal={Duke Math. J.},
   volume={53},
   date={1986},
   pages={869--943},
}
\bib{Ouh}{book}{
   author={Ouhabaz, El Maati},
   title={Analysis of heat equations on domains},
   series={London Mathematical Society Monographs Series},
   volume={31},
   publisher={Princeton University Press, Princeton, NJ},
   date={2005},
   pages={xiv+284},
}
\bib{P3}{article}{
   author={Pinchover, Yehuda},
   title={On criticality and ground states of second order elliptic
   equations. II},
   journal={J. Differential Equations},
   volume={87},
   date={1990},
   pages={353--364},
}
\bib{P0}{article}{
   author={Pinchover, Yehuda},
   title={Large time behavior of the heat kernel and the behavior of the
   Green function near criticality for nonsymmetric elliptic operators},
   journal={J. Funct. Anal.},
   volume={104},
   date={1992},
   pages={54--70},
}
\bib{P1}{article}{
   author={Pinchover, Yehuda},
   title={On positivity, criticality, and the spectral radius of the shuttle
   operator for elliptic operators},
   journal={Duke Math. J.},
   volume={85},
   date={1996},
   pages={431--445},
}
\bib{P2}{article}{
   author={Pinchover, Yehuda},
   title={Large time behavior of the heat kernel},
   journal={J. Funct. Anal.},
   volume={206},
   date={2004},
   pages={191--209},
}
\bib{PZ}{article}{
   author={Pinchover, Yehuda},
   title={Some aspects of large time behavior of the heat kernel: an
   overview with perspectives},
   conference={
      title={Mathematical physics, spectral theory and stochastic analysis},
   },
   book={
      series={Oper. Theory Adv. Appl.},
      volume={232},
      publisher={Birkh\"{a}user/Springer Basel AG, Basel},
   },
   date={2013},
   pages={299--339},
}
\bib{S}{article}{
   author={Simon, Barry},
   title={Large time behavior of the $L^{p}$ norm of Schr\"{o}dinger
   semigroups},
   journal={J. Functional Analysis},
   volume={40},
   date={1981},
   pages={66--83},
}
\bib{Zhang0}{article}{
   author={Zhang, Qi S.},
   title={Large time behavior of Schr\"{o}dinger heat kernels and applications},
   journal={Comm. Math. Phys.},
   volume={210},
   date={2000},
   pages={371--398},
}
\bib{Zhang}{article}{
   author={Zhang, Qi S.},
   title={Global bounds of Schr\"{o}dinger heat kernels with negative
   potentials},
   journal={J. Funct. Anal.},
   volume={182},
   date={2001},
   pages={344--370},
}
\end{biblist}
\end{bibdiv}
\end{document}